\documentclass{amsart}
\usepackage{amsmath,amssymb,amsthm,mathtools}
\usepackage[hidelinks]{hyperref}

\theoremstyle{plain}
\newtheorem{theorem}{Theorem}[section]
\newtheorem{lemma}[theorem]{Lemma}
\newtheorem{proposition}[theorem]{Proposition}
\newtheorem{corollary}[theorem]{Corollary}
\theoremstyle{definition}
\newtheorem{remark}[theorem]{Remark}
\newtheorem{example}[theorem]{Example}

\theoremstyle{plain}
\newtheorem{theoremA}{Theorem}

\newcommand{\intZ}{\mathbb{Z}}
\newcommand{\ratQ}{\mathbb{Q}}
\newcommand{\ratQp}{\ratQ_p}
\newcommand{\intring}[1]{O_{#1}}
\newcommand{\OK}{\intring{K_b}}

\newcommand{\disc}[1]{D_{#1}}
\DeclareMathOperator{\Norm}{N}
\DeclareMathOperator{\Gm}{\mathbf{G}_m}
\DeclareMathOperator{\Jac}{Jac}

\newcommand{\Dt}{\widetilde{D}}
\newcommand{\Ct}{\widetilde{C}}
\newcommand{\Pt}{\widetilde{P}}
\newcommand{\Qbar}{\overline{\ratQ}}

\begin{document}

\title[Explicit units and curves from their powers]
{On explicit units in a family of cubic fields and the arithmetic of curves
arising from powers of these units}

\author{Iwao Kimura}
\address{Department of Mathematics, University of Toyama,
3190 Gofuku, Toyama 930-8555, Japan}
\email{iwao@sci.u-toyama.ac.jp}

\author{Hikaru Umemoto}
\address{Graduate School of Science and Engineering, University of Toyama,
3190 Gofuku, Toyama 930-8555, Japan}
\email{m25c1008@ems.u-toyama.ac.jp}


\begin{abstract}
For an integer $b\neq0,1$ let $K_b$ be the cubic field generated by the unique
real root of $x^3-3bx-b^3$, let $\eta_0$ be its fundamental unit, and let
$\eta_b$ be the explicit unit studied in the companion paper. Write
$\eta_b=\eta_0^{\,n(b)}$. Deciding whether $n(b)$ is divisible by $q$ amounts to
finding the integral points of a plane curve $C_q$ of degree $2q$. We prove that
for \emph{every} integer $q\ge1$ the curve $C_q$ is irreducible with
\[
  g(C_q)=\binom{q-1}{2},
\]
the genus of a smooth plane curve of degree $q$ rather than of degree $2q$. The
proof is not a genus computation but a structural one: $C_q$ is the inverse
image, under the $q$-power map of a two-dimensional torus, of the curve traced
by the minimal polynomial of $\eta_b$ as $b$ varies, and the
formula falls out of two applications of Riemann--Hurwitz in the resulting
tower. As a consequence the problem splits into exactly three regimes: $q=2$
rational, $q=3$ elliptic, and $q\ge4$ of general type. In particular only
finitely many $b$ have $q\mid n(b)$ for each $q\ge4$. A lower bound for linear
forms in two logarithms gives $n(b)<6082138$ for every $b$, so $n(b)\in\{1,2\}$
for all but finitely many $b$, though not effectively. We add two computations:
$n(b)$ is determined outright for $|b|\le10^{30}$, and the rational points of
$C_4$ are determined by quadratic Chabauty, giving
$4\mid n(b)\iff b\in\{-75,-3,5,96\}$. We state precisely which parts of the
latter rest on certified arithmetic and which do not. Finally, we show that
$K_b\cong K_{b'}$ with $b\neq b'$ holds for $\{b,b'\}=\{-1,5\}$, and otherwise
only if one of the two parameters, $c$ say, satisfies $|c|>10^{30}$ and $n(c)$
is divisible by a prime $q\ge5$. There are only finitely many such pairs.
\end{abstract}

\maketitle

\bibliographystyle{amsplain-doi}

\section{Introduction}\label{sec:intro}

\subsection{The family}\label{sec:family}
Let $\intZ$ be the ring of rational integers, and $\ratQ\subset\mathbb{R}$ the
fields of rational and of real numbers. For $b\in\intZ$ put $f_b(x)=x^3-3bx-b^3$, of discriminant
$-27b^3(b^3-4)$; this is negative exactly for $b\neq0,1$, and \emph{throughout,
$b$ denotes an integer with $b\neq0,1$}. Let $\theta$ be the unique real root of
$f_b$ and $K_b=\ratQ(\theta)$, a cubic field with one real and one pair of complex
embeddings, so that the unit group $\OK^\times$ of its ring of integers has rank
one. Let $\eta_0>1$ be the fundamental unit, and set
\[
  \eta_b\ :=\ \frac{1}{(b+1)-\theta} .
\]

The family and the numerical observations behind it come from the companion
paper \cite{KimuraUmemoto}. It shows that $\eta_b$ is a unit greater than $1$;
that it is fundamental for an explicit set of $b$; and that $\eta_b$ is a square
in $\OK^\times$ exactly for the $b$ of a family parametrized by the Pell equation
$D^2-3E^2=1$, which is therefore infinite. Under an explicit condition on the
discriminant of $K_b$ these square-exceptions are all the $b$ for which $\eta_b$
fails to be fundamental.

Write
\[
  \eta_b=\eta_0^{\,n(b)},\qquad n(b)\ge1
\]
(Lemma~\ref{lem:eta}). In these terms, \cite{KimuraUmemoto} shows that $n(b)=1$
for an explicit set of $b$, and that $n(b)$ is even exactly for the $b$ of the
Pell family. Under the condition on the discriminant, $n(b)=1$ for every other
$b$. It leaves open the value of $n(b)$ when $n(b)$ is even, and the case in
which the condition fails. The present paper is about the rest: for which $b$, and which $q$, is
$q\mid n(b)$?

\subsection{The curves}\label{sec:curves}
The minimal polynomial of $\eta_b$ is
\begin{equation}\label{eq:etamin}
  X^3-3(b^2+b+1)X^2+3(b+1)X-1 .
\end{equation}
A unit $u>1$ of $K_b$ has characteristic polynomial $X^3-BX^2+AX-1$ with
$A,B\in\intZ$ (\S\ref{sec:setup}), so the condition $\eta_b=u^q$ becomes a pair of
equations in $(A,B)$, and eliminating $b$ leaves one. Define $p_q\in\intZ[A,B]$ by
Newton's identities,
\begin{equation}\label{eq:newton}
  p_0=3,\quad p_1=B,\quad p_2=B^2-2A,\qquad
  p_{k+1}=B\,p_k-A\,p_{k-1}+p_{k-2}\quad(k\ge2),
\end{equation}
Then $p_q(A,B)=\sum_ie_i^{\,q}$ and $p_q(B,A)=\sum_{i<j}(e_ie_j)^q$ for the
roots $e_1,e_2,e_3$ of $X^3-BX^2+AX-1$. Put
\begin{equation}\label{eq:Cq}
  C_q\ :\quad F_q(A,B):=p_q(B,A)^2-3\,p_q(B,A)+9-3\,p_q(A,B)=0 .
\end{equation}
Then $q\mid n(b)$ if and only if $C_q$ carries a suitable integral point;
Proposition~\ref{prop:dictionary} makes this precise.

\subsection{Results}\label{sec:results}
Our main theorem is about the geometry of $C_q$. Here $q$ is any positive
integer, prime or not.

\begin{theoremA}\label{thmA}
For every $q\ge1$ the polynomial $F_q$ is squarefree, the curve $C_q$ is
geometrically irreducible of degree $2q$, and its smooth projective model has
genus
\[
  g(C_q)=\binom{q-1}{2} .
\]
\end{theoremA}

This is the genus of a smooth plane curve of degree $q$, whereas $C_q$ has degree
$2q$ and arithmetic genus $(2q-1)(q-1)$; the deficit
$\sum_P\delta_P=\tfrac32q(q-1)$ is thus forced. For $q=4$ the curve is indeed a
smooth plane quartic (\S\ref{sec:c4curve}), but that does not explain the formula,
and neither does the Fermat curve of degree $q$, which has the same genus but is
not $C_4$ (Remark~\ref{rem:toric}). The explanation is instead the following.
We write $\mathbf{A}^n$ and $\mathbf{P}^n$ for affine and projective $n$-space,
with the coordinates or the base ring as a subscript, and $\Gm$ for the
multiplicative group. Let $T$ be the torus $\{e_1e_2e_3=1\}$ and identify $\mathbf{A}^2_{A,B}$ with $T/S_3$
by sending a triple to the coefficients of $\prod_i(X-e_i)$; then
$(A,B)\mapsto(p_q(B,A),p_q(A,B))$ is the $S_3$-quotient of the $q$-power map of
$T$, and $C_q$ is the inverse image of the curve $D$ traced by the coefficients
of \eqref{eq:etamin} as $b$ varies. Riemann--Hurwitz applied twice
in the resulting tower yields the formula for all $q$ at once. The details occupy
\S\S\ref{sec:torus}--\ref{sec:descent}.

Theorem~\ref{thmA} separates the Diophantine problem into three regimes.

\begin{theoremA}\label{thmB}
\begin{enumerate}
\item $g(C_2)=0$, and there are infinitely many $b$ with $2\mid n(b)$.
\item $g(C_3)=1$, and $b=9$ is the only integer with $3\mid n(b)$.
\item $g(C_q)\ge2$ for $q\ge4$, so only finitely many $b$ satisfy $q\mid n(b)$.
\item $n(b)<6082138$ for every $b\neq0,1$. Hence $n(b)\in\{1,2\}$ for all but
  finitely many $b$.
\end{enumerate}
\end{theoremA}

Parts (i) and (ii) are results of \cite{KimuraUmemoto}, proved there by a Pell
equation and by a Thue equation respectively; Theorem~\ref{thmA} accounts for the
shape of both and supplies (iii). Part (iii) is not effective.
Theorem~\ref{thmA} makes it uniform in $q$: previously the genus was available
only for those $q$ for which it had been computed. Part (iv) combines (ii) and
(iii) with a bound for $n(b)$ that does not depend on $b$. We obtain this bound
from a lower bound for linear forms in two logarithms, the tool we also use for
Theorem~\ref{thmE}. The bound is effective. The finiteness in (iv) is not, since
it rests on (iii).

The remaining three results rest, at least in part, on computation, and we
describe their status in Remark~\ref{rem:status}.

\begin{theoremA}\label{thmC}
Let $b\neq0,1$ with $|b|\le10^{30}$. Then $n(b)\in\{1,2,4,6,8,13\}$; moreover
$n(b)=4$ exactly for $b=-75,5,96$, $n(b)=6$ exactly for $b=9$, $n(b)=8$ exactly
for $b=-3$, $n(b)=13$ exactly for $b=3$, and $n(b)=2$ exactly for the remaining
square-exceptions.
\end{theoremA}

\begin{theoremA}\label{thmD}
For every $b\neq0,1$, $\;4\mid n(b)$ if and only if $b\in\{-75,-3,5,96\}$.
\end{theoremA}

Finally, the map $b\mapsto K_b$ is not injective, but it is close to being so.

\begin{theoremA}\label{thmE}
\begin{enumerate}
\item $K_{-1}=K_5$; indeed $\eta_5=\eta_{-1}^{\,4}$.
\item Let $m\ge2$. The minimal polynomial of $\eta_b^{\,m}$ equals that of some
  $\eta_c$ only for $(b,m,c)=(-1,4,5)$. Consequently, if $n(b)=1$ and $b\neq-1$,
  then $K_b\not\cong K_{b'}$ for every $b'\neq b$; and $K_{-1}=K_5$ is $K_c$ for no
  other $c$.
\item Let $b\neq b'$ with $K_b\cong K_{b'}$ and $\{b,b'\}\neq\{-1,5\}$. Then
  for some $c\in\{b,b'\}$ we have $|c|>10^{30}$, and $n(c)$ is divisible by a
  prime $q\ge5$. There are only finitely many such pairs $\{b,b'\}$.
\end{enumerate}
\end{theoremA}

\begin{remark}\label{rem:status}
Theorem~\ref{thmA} and Theorem~\ref{thmB}(iii) have proofs that use no machine
computation. Parts (i) and (ii) of Theorem~\ref{thmB} are quoted from
\cite{KimuraUmemoto}, where part (ii) rests on the solution of a cubic Thue
equation. Part (iv) rests in addition on Smyth's lower bound for the Mahler
measure of a non-reciprocal algebraic integer (Lemma~\ref{lem:smyth}), on a
theorem of Laurent, Mignotte and Nesterenko \cite{LMN1995} on linear forms in
two logarithms, and on one inequality that we check in exact rational
arithmetic (Proposition~\ref{iso:uniform}). Theorem~\ref{thmC} rests on a
finite computation in exact integer arithmetic (\S\ref{sec:sweep}). Smyth's
lower bound also limits the exponents to be examined there.
Theorem~\ref{thmD} rests on a quadratic Chabauty computation for $C_4$ at
$p=13$ and on a Mordell--Weil sieve. Both are exact: the arithmetic is carried
out modulo powers of $p$ and in finite fields, every stage carries a proved
precision, and no zero is lost beyond that precision. We construct and check
here the two regular models they use, at the primes $2$ and $3$; the algorithm
we use at $2$ is unpublished. See \S\ref{rem:certificate} and
\S\ref{sec:precision}. Theorem~\ref{thmD} also follows from an independent
argument by elliptic logarithms (Remark~\ref{rem:elllog}), which rests on
theorems of David and of Silverman, on the Mordell--Weil group of $E$ in
\eqref{eq:E} and on exact arithmetic, and uses neither quadratic Chabauty nor
the regular models.
Theorem~\ref{thmE}(i) is proved by hand. Part~(ii) rests on the theorem of
Laurent, Mignotte and Nesterenko \cite{LMN1995} mentioned above, and on a
finite computation in exact rational arithmetic
(Remark~\ref{iso:status}); its last assertion also uses $n(-1)=1$, from
\cite{KimuraUmemoto}. The first assertion of part~(iii) rests in addition on
Theorems~\ref{thmB}(ii), \ref{thmC} and~\ref{thmD}, and on the computation of
Lemma~\ref{iso:six}. The finiteness of the pairs rests instead on part~(ii) and
Theorem~\ref{thmB}(ii)--(iv). It does not use Theorems~\ref{thmC} or~\ref{thmD},
and it is not effective.
\end{remark}

\subsection{Relation to known constructions}\label{sec:known}
The map $\Phi_q$ of \S\ref{sec:torus} is the Adams operation $\psi^q$ read on the
coefficients of a characteristic polynomial, and the identification
$\mathbf{A}^2\cong T/S_3$ is the adjoint quotient of $\mathrm{SL}_3$
\cite[\S6]{Steinberg1965}. In that setting it is standard that the $q$-power map
of a torus is \'etale of degree $q^{\operatorname{rank}}$ in characteristic
zero. The
substitution itself, which replaces the root of a defining polynomial by its
$q$-th power so as to introduce a second parameter, occurs in
L\'evesque--Waldschmidt \cite{LevesqueWaldschmidt2015} for Shanks's simplest
cubic fields, where the resulting family of Thue equations is solved effectively;
the geometry of the associated curves is not considered there. Here we prove
the genus formula, and for that the standard toric machinery is unavailable: the
lattice-point bound of \cite{Khovanskii1977,Khovanskii1978} is attained
neither by $\Dt$ nor by $C_q$ for $q\ge2$ with $3\nmid q$ (Remark~\ref{rem:toric}).

\subsection{Outline}\label{sec:outline}
\S\ref{sec:setup} fixes the correspondence between units and integral points.
\S\ref{sec:torus} sets up the torus and the tower, \S\ref{sec:base} treats its
base curve, \S\ref{sec:kummer} the Kummer cover, and \S\ref{sec:descent} the
descent along $S_3$, which completes the proof of Theorem~\ref{thmA}.
Theorem~\ref{thmB}(i)--(iii) and Theorems~\ref{thmC} and~\ref{thmD} are proved in
\S\S\ref{sec:consequences}--\ref{sec:c4}, and Theorem~\ref{thmB}(iv) and
Theorem~\ref{thmE} in \S\ref{sec:iso}.

\section{Setup}\label{sec:setup}

We collect what is needed from the companion paper. Recall that $\OK$ denotes the
ring of integers of $K_b$ and that $\OK^\times$ has rank one, because $K_b$ has
one real and one pair of complex embeddings. The fundamental unit is $\eta_0>1$,
and $\Norm$ is the norm of $K_b/\ratQ$.

\begin{lemma}\label{lem:eta}
Let $b$ be an integer, $b\neq0,1$. Then $\eta_b$ is a unit of $\OK$ with
$\eta_b>1$ and $\Norm(\eta_b)=1$, its minimal polynomial is \eqref{eq:etamin},
and $\eta_b=\eta_0^{\,n(b)}$ with $n(b)\ge1$.
\end{lemma}

\begin{proof}
Eliminating $\theta$ between $f_b(\theta)=0$ and $\eta_b\bigl((b+1)-\theta\bigr)=1$
gives \eqref{eq:etamin}, which is monic with constant term $-1$; so $\eta_b$ is a
unit and $\Norm(\eta_b)=1$. As $f_b(b)=-3b^2<0$ and $f_b(b+1)=1>0$, the unique
real root satisfies $b<\theta<b+1$, so $0<(b+1)-\theta<1$ and $\eta_b>1$. Finally
$\OK^\times=\{\pm1\}\times\langle\eta_0\rangle$, and $\eta_b>1$ leaves only
$\eta_b=\eta_0^{\,k}$ with $k\ge1$.
\end{proof}

For a unit $u>1$ of $K_b$ one has $\Norm(u)=u\,|\sigma(u)|^2>0$, where $\sigma$ is
either complex embedding, so $\Norm(u)=1$; and $u\notin\ratQ$, since the only
rational units are $\pm1$. Hence $u$ generates $K_b$ and its characteristic
polynomial is $X^3-BX^2+AX-1$ with $A,B\in\intZ$. This is the source of the
coordinates $(A,B)$ used throughout.

\begin{proposition}\label{prop:dictionary}
Let $b$ be an integer, $b\neq0,1$, and let $q\ge1$. Then $q\mid n(b)$ if and only
if there is a point $(A,B)\in\intZ^2$ on $C_q$ with
\begin{equation}\label{eq:fibre}
  p_q(B,A)=3(b+1) .
\end{equation}
Moreover the point is unique: $X^3-BX^2+AX-1$ is the characteristic polynomial of
$\eta_0^{\,n(b)/q}$.
\end{proposition}

\begin{proof}
Suppose $q\mid n(b)$ and put $u=\eta_0^{\,n(b)/q}$, so that $u>1$ and
$\eta_b=u^q$. Let $e_1,e_2,e_3$ be the conjugates of $u$ and $X^3-BX^2+AX-1$ its
characteristic polynomial, with $A,B\in\intZ$ by the paragraph before
Proposition~\ref{prop:dictionary}. The
$e_i^{\,q}$ are the conjugates of $\eta_b$, so comparing with \eqref{eq:etamin},
\begin{equation}\label{eq:pair}
  p_q(A,B)=\sum_i e_i^{\,q}=3(b^2+b+1),\qquad
  p_q(B,A)=\sum_{i<j}(e_ie_j)^q=3(b+1) .
\end{equation}
Substituting \eqref{eq:pair} in \eqref{eq:D} gives $F_q(A,B)=0$.

Conversely, let $(A,B)\in\intZ^2$ satisfy $F_q(A,B)=0$ and \eqref{eq:fibre}, and
let $e_1,e_2,e_3$ be the roots of $X^3-BX^2+AX-1$, algebraic integers of product
$1$ and hence units. The second identity of \eqref{eq:pair} holds by hypothesis
and the first because $F_q(A,B)=0$. Hence the $e_i^{\,q}$ are the roots of
\eqref{eq:etamin}, that is, the conjugates of $\eta_b$. Two of these are not
real, so neither are two of the $e_i$; as $X^3-BX^2+AX-1$ has real coefficients,
exactly one root is real. Let $e_1$ be the real one, so that $e_1^{\,q}=\eta_b$.
Since $e_2,e_3$ are complex conjugate, $e_1=|e_2|^{-2}>0$, so $e_1$ is the positive
$q$-th root of $\eta_b$. The field $\ratQ(e_1)$ has degree at most $3$ and contains
$\eta_b$, which generates $K_b$; so $\ratQ(e_1)=K_b$, and $u:=e_1$ is a unit
$u>1$ of $K_b$ with $u^q=\eta_b$, whence $q\mid n(b)$, and $u=\eta_0^{\,n(b)/q}$
determines $(A,B)$.
\end{proof}

Condition \eqref{eq:fibre} says that $(A,B)$ lies in the fiber of
$C_q\to\mathbf{A}^1_b$ over $b$; it is not implied by $F_q(A,B)=0$, which only
places $(A,B)$ over \emph{some} $b$. Thus Proposition~\ref{prop:dictionary}
converts a question about units into one about integral points, and
Theorem~\ref{thmA} controls the geometry of the curve on which they lie. For
$q=1$ it degenerates: $p_1(B,A)=A$, so $C_1=D$ and \eqref{eq:fibre} reads
$A=3(b+1)$.

\section{The torus and the \texorpdfstring{$q$}{q}-power map}\label{sec:torus}

\subsection*{Conventions}
In the proof of Theorem~\ref{thmA}, \S\S\ref{sec:torus}--\ref{sec:descent}, we work
with curves over $\Qbar$, a fixed algebraic closure of $\ratQ$, and $\mu_n\subset\Qbar^\times$ is the group of $n$-th roots of unity.
The letter $p$ always denotes a prime number. For a smooth
projective curve $Y$ over $\Qbar$ we write $g(Y)$ for its genus,
$\chi(Y)=2-2g(Y)$ for its Euler characteristic, $\Qbar(Y)$ for its function
field, and $\operatorname{ord}_P$ for the normalized valuation at a point $P$.
For an open subset $Y\setminus S$ with $S$ finite we put
$\chi(Y\setminus S)=\chi(Y)-\#S$.
A curve given by an affine equation is identified with its smooth projective
model, and a superscript $\circ$ marks an open subset of it.

Let
\[
  T=\{(e_1,e_2,e_3)\in\Gm^3\ :\ e_1e_2e_3=1\}
\]
be the two-dimensional split torus over $\ratQ$, with
\[
  X^*(T)=\intZ^3/\intZ(1,1,1),\qquad
  X_*(T)=\{n\in\intZ^3:n_1+n_2+n_3=0\}
\]
as character and cocharacter lattices, the pairing $\langle\ ,\ \rangle$
being induced by the standard pairing on $\intZ^3$. A
character $\chi=(a_1,a_2,a_3)$ is the function
$e\mapsto e_1^{a_1}e_2^{a_2}e_3^{a_3}$ on $T$, well defined because
$e_1e_2e_3=1$.

The symmetric group $S_3$ acts on $T$ by permuting the coordinates, and
\begin{equation}\label{eq:pi}
  \pi\colon T\longrightarrow \mathbf{A}^2_{A,B},\qquad
  \pi(e)=(A,B)=\Bigl(\sum_{i<j}e_ie_j,\ \sum_i e_i\Bigr),
\end{equation}
identifies $\mathbf{A}^2_{A,B}$ with the quotient $T/S_3$; thus $\pi$ is finite of
degree $6$, and the fiber of $\pi$ through $e$ consists of the roots of
$X^3-BX^2+AX-1$. This is the adjoint quotient of $\mathrm{SL}_3$
\cite[\S6]{Steinberg1965}.

Let $[q]\colon T\to T$ be $e\mapsto(e_1^q,e_2^q,e_3^q)$. It is an isogeny of
degree $q^2$, and is \'etale because the characteristic is zero. It commutes
with the action of $S_3$, hence descends to a finite map
$\Phi_q\colon\mathbf{A}^2\to\mathbf{A}^2$ with $\Phi_q\circ\pi=\pi\circ[q]$, and
$\deg\Phi_q=q^2$. In the coordinates \eqref{eq:pi} one has
$\Phi_q(A,B)=(p_q(B,A),p_q(A,B))$ with $p_q$ as in \eqref{eq:newton}.

Let $D\subset\mathbf{A}^2$ be the image of the map $b\mapsto(3(b+1),3(b^2+b+1))$
attached to \eqref{eq:etamin}. Eliminating $b$ shows that $D$ is the plane curve
\begin{equation}\label{eq:D}
  A^2-3A+9-3B=0 ,
\end{equation}
and $b\mapsto A$ is an isomorphism $\mathbf{A}^1_b\xrightarrow{\ \sim\ }D$. By
\eqref{eq:Cq},
\begin{equation}\label{eq:CqPhi}
  C_q=\Phi_q^{-1}(D) .
\end{equation}
Put
\[
  \Dt=\pi^{-1}(D)\subset T,\qquad \Ct_q=[q]^{-1}(\Dt)\subset T ,
\]
so that $\Dt$ is the curve of ordered triples of conjugates of $\eta_b$. Write
$\Dt^\circ=\Dt\cap T$ and $\Ct_q^\circ=\Ct_q\cap T=[q]^{-1}(\Dt^\circ)$ for the
parts of these curves lying in the torus. Finally, these curves form the
commutative diagram
\begin{equation}\label{eq:tower}
  \begin{array}{ccc}
    \Ct_q & \xrightarrow{\ [q]\ } & \Dt\\[2pt]
    \downarrow & & \downarrow\\[2pt]
    C_q & \xrightarrow{\ \Phi_q\ } & D=\mathbf{A}^1_b ,
  \end{array}
\end{equation}
in which the vertical maps are the quotients by $S_3$.

\section{The base curve \texorpdfstring{$\Dt$}{D-tilde}}\label{sec:base}

\begin{proposition}\label{prop:base}
\begin{enumerate}
\item $\Dt$ is geometrically irreducible of genus $1$, and $\Dt\to\mathbf{P}^1_b$
  is Galois with group $S_3$.
\item Its branch locus is $\{0\}\cup\{b:b^3=4\}$. Over each of these four points
  lie three points of $\Dt$, of ramification index $2$, whose stabilizers are the
  three transpositions.
\item The fiber over $b=\infty$ consists of six points, on which $S_3$ acts
  simply transitively; in particular no nontrivial element of $S_3$ fixes one of
  them. At each such point $P$ the vector
  $v_P:=\bigl(\operatorname{ord}_Pe_1,\operatorname{ord}_Pe_2,\operatorname{ord}_Pe_3\bigr)$
  lies in $X_*(T)$ and is a permutation of $(-2,1,1)$.
\item The functions $e_1,e_2,e_3$ have neither zeros nor poles outside those six
  points. Consequently $\Dt^\circ$ is the complement of the six, and
  $\chi(\Dt^\circ)=-6$.
\end{enumerate}
\end{proposition}

\begin{proof}
Write $\Dt_1=\Dt/\langle(2\,3)\rangle$, the curve of pairs $(x,b)$ with $x$ a root
of \eqref{eq:etamin}; it has degree $3$ over $\mathbf{A}^1_b$. Read as an
equation in $b$, \eqref{eq:etamin} is
$3x^2b^2+3x(x-1)b-(x-1)^3=0$, of discriminant $3x^2(x-1)^2(4x-1)$. Hence $\Dt_1$
is the conic $s^2=3(4x-1)$, parametrized by
\begin{equation}\label{eq:param}
  x=\frac{s^2+3}{12},\qquad b=\frac{(s-3)^2(s+3)}{6(s^2+3)} ,
\end{equation}
so $\Dt_1\cong\mathbf{P}^1_s$ and $g(\Dt_1)=0$.

The discriminant of \eqref{eq:etamin} with respect to $X$ is $-27b^3(b^3-4)$,
which gives the branch locus in (ii). Substituting $X=1+u$ in \eqref{eq:etamin}
produces $u^3-3bu-3b^2$ plus terms lying above the Newton polygon of these
three; the polygon has slopes $\tfrac12,\tfrac12,1$, so two of the roots are
$u=\pm\sqrt{3b}+O(b)$ and the third is $u=-b+O(b^2)$. The monodromy at $b=0$
is therefore a transposition. As $b^3-4$ has three simple zeros, which are
simple zeros of the discriminant, the monodromy at each of them is a
transposition as well.

For $b=\infty$ put $t=1/b$. The Newton polygon of $t^2f(Y,1/t)$, where $f$
denotes \eqref{eq:etamin}, has slopes $-2$ and $+1$, the latter twice; and after
$Y=tZ$ the slope $+1$ part reduces at $t=0$ to $-3Z^2+3Z-1$, whose discriminant
$-3$ is nonzero. The two roots of order $+1$ thus stay distinct, and separating
them requires only the constant extension $\ratQ(\sqrt{-3})$. Hence $b=\infty$
is unramified, and the orders of $e_1,e_2,e_3$ at a point above it are $-2,1,1$
in some order.

Since $b^3$ has odd order at $b=0$, the discriminant is not a square in
$\Qbar(b)$; as $\Dt_1$ is irreducible, the Galois group of \eqref{eq:etamin} over
$\Qbar(b)$ is $S_3$, which is (i) apart from the genus. Being Galois of degree
$6$ and unramified at $b=\infty$, the cover has six points there, permuted simply
transitively; this is (iii), and it gives (ii) because the stabilizer of a point
over a branch point is generated by the local monodromy. Riemann--Hurwitz now
reads $2g(\Dt)-2=6\cdot(-2)+12$, so $g(\Dt)=1$.

For (iv), the $e_i$ are roots of the monic polynomial \eqref{eq:etamin}, whose
coefficients lie in $\Qbar[b]$ and whose constant term is a unit; so the $e_i$
are units in the integral closure of $\Qbar[b]$ in $\Qbar(\Dt)$ and have no zeros
or poles over finite $b$. Removing six points from a curve of genus $1$ leaves
Euler characteristic $-6$.
\end{proof}

\section{The Kummer cover \texorpdfstring{$\Ct_q\to\Dt$}{C-tilde\_q -> D-tilde}}
\label{sec:kummer}

Over $\Dt^\circ$ the cover $\Ct_q\to\Dt$ is obtained by adjoining $q$-th roots of
$e_1$ and $e_2$; it is the $(\intZ/q)^2$-Kummer cover attached to the
homomorphism
\begin{equation}\label{eq:kummer}
  \kappa_q\colon X^*(T)/q\longrightarrow L^\times/(L^\times)^q,
  \qquad L:=\Qbar(\Dt),
\end{equation}
which sends a character to the corresponding function. Since $\Qbar$ contains
$\mu_q$, Kummer theory \cite[Ch.~VI, \S8]{Lang2002} identifies the Galois group
of $L(e_1^{1/q},e_2^{1/q})/L$ with the dual of the image of $\kappa_q$; so the
cover is connected if and only if $\kappa_q$ is injective, and we show this
first.

\begin{lemma}\label{lem:cube}
$e_2/e_3$ is not a cube in $L^\times$.
\end{lemma}

\begin{proof}
Put $u=e_2/e_3$ and suppose $u=h^3$ with $h\in L^\times$.

By Proposition~\ref{prop:base}(iii), (iv), $\operatorname{ord}_P u$ takes the values
$0,0,-3,-3,3,3$ at the six points over $b=\infty$ and vanishes elsewhere. Hence
$\operatorname{div}_\infty u$ has degree $6$, and $\operatorname{div}h=\tfrac13\operatorname{div}u$
shows that $h$ has degree $2$ as a map to $\mathbf{P}^1$.

Let $\iota=(2\,3)\in S_3$. Then $\iota u=u^{-1}$, so $(\iota(h)\,h)^3=1$; a
function whose cube is $1$ is constant, and we may write $\iota(h)\,h=\zeta$ with
$\zeta\in\mu_3$.

The quotient map $\Dt\to\Dt/\iota=\Dt_1$ has degree $2$, and
$g(\Dt)=1$, $g(\Dt_1)=0$ by Proposition~\ref{prop:base} and its proof. By
Riemann--Hurwitz $\iota$ has exactly four fixed points $F_1,\dots,F_4$. None of
them lies over $b=\infty$, by Proposition~\ref{prop:base}(iii); and
$\operatorname{div}h=\tfrac13\operatorname{div}u$ is supported over $b=\infty$ by
(iv) of the same proposition, so $h(F_i)\in\Qbar^\times$.

Since $\iota$ fixes $F_i$ and interchanges $e_2$ and $e_3$, we have
$e_2(F_i)=e_3(F_i)$ and therefore $u(F_i)=1$. Thus $h(F_i)^3=1$, while
$\iota(h)h=\zeta$ evaluated at $F_i$ gives $h(F_i)^2=\zeta$. Dividing,
$h(F_i)=\zeta^{-1}$ for every $i$. A map of degree $2$ takes no value at more
than two points, whereas $h$ takes the value $\zeta^{-1}$ at four.
\end{proof}

\begin{lemma}\label{lem:kernel}
$\kappa_q$ is injective for every $q\ge1$.
\end{lemma}

\begin{proof}
Let $\chi\in X^*(T)$ represent an element of $\ker\kappa_q$. Then
$\operatorname{div}(\chi)\in q\operatorname{Div}(\Dt)$, so
$\langle\chi,v_P\rangle\equiv0\pmod q$ at each of the six points over
$b=\infty$. By Proposition~\ref{prop:base}(iii) the $v_P$ run through the
permutations of $(-2,1,1)$, and these span a sublattice of $X_*(T)$ of index $3$.
Hence $\ker\kappa_q$ is trivial when $3\nmid q$.

Let $q=3m$. The annihilator of that sublattice in $X^*(T)/q$ is then cyclic of
order $3$, generated by the class of $m\chi_1$ with $\chi_1=(1,2,0)$; and
$\chi_1$ is the function $e_1e_2^2=e_2/e_3$. Now
$(e_2/e_3)^m\in(L^\times)^q$ forces $e_2/e_3=\zeta h^3$ for some $h\in L^\times$
and some $\zeta\in\mu_m$, and $\zeta$ is a cube in $\Qbar^\times$; so $e_2/e_3$
would be a cube in $L^\times$, contradicting Lemma~\ref{lem:cube}.
\end{proof}

\begin{proposition}\label{prop:kummer}
For every $q\ge1$ the curve $\Ct_q$ is geometrically irreducible, the six points
of $\Dt$ over $b=\infty$ each carry $q$ points of $\Ct_q$ with ramification index
$q$, and
\[
  g(\Ct_q)=3q^2-3q+1 .
\]
\end{proposition}

\begin{proof}
As $[q]$ is \'etale, so is $\Ct_q^\circ\to\Dt^\circ$, of degree $q^2$; therefore
$\chi(\Ct_q^\circ)=q^2\chi(\Dt^\circ)=-6q^2$ by
Proposition~\ref{prop:base}(iv). Let $P$ lie over $b=\infty$, with local
parameter $\tau$, and write $e_i=\tau^{n_i}w_i$ with $w_i$ a unit in the
completed local ring. The residue field $\Qbar$ is algebraically closed of
characteristic $0$, so by Hensel's lemma each $w_i$ is a $q$-th power in the
completed local ring, and the completion of $\Ct_q$ at a point above $P$ is
$\Qbar((\tau))\bigl[\tau^{n_1/q},\tau^{n_2/q}\bigr]$; this is totally ramified of
degree $q/\gcd(q,n_1,n_2,n_3)$ \cite[\S3.7]{Stichtenoth2009}. By
Proposition~\ref{prop:base}(iii) the vector $v_P=(n_1,n_2,n_3)$ is a permutation of
$(-2,1,1)$ and so is primitive in $X_*(T)$, whence the degree is $q$: there are
$q$ points above $P$, each with ramification index $q$. Adding the resulting $6q$
points,
\[
  \chi(\Ct_q)=-6q^2+6q .
\]
By Lemma~\ref{lem:kernel} the cover is connected, so $\Ct_q$ is irreducible and
$2-2g(\Ct_q)=-6q^2+6q$.
\end{proof}

\section{The \texorpdfstring{$S_3$}{S\_3} descent}\label{sec:descent}

\begin{proposition}\label{prop:descent}
$F_q$ is squarefree, $C_q$ is geometrically irreducible of degree $2q$, the map
$\Ct_q\to C_q$ has degree $6$, and it is ramified at exactly $12q$ points, each
with ramification index $2$.
\end{proposition}

\begin{proof}
If a point $t\in T$ has $[q](t)=e$ with $e_1,e_2,e_3$ distinct, then the
coordinates of $t$ are distinct as well, so $\pi$ is \'etale at $t$ and at $e$.
As $[q]$ is \'etale, $\Phi_q\circ\pi=\pi\circ[q]$ is \'etale at $t$; and a
morphism $g$ is \'etale at $f(x)$ whenever $f$ and $g\circ f$ are \'etale at $x$,
so $\Phi_q$ is \'etale at $\pi(t)$. Since $\Phi_q$ is finite,
every component of $C_q=\Phi_q^{-1}(D)$ dominates $D$, and its generic point lies
over the generic point of $D$, where the $e_i$ are distinct. Hence $\Phi_q^*(D)$
is reduced along every component and $F_q$ is squarefree. By
Proposition~\ref{prop:kummer} and \eqref{eq:tower}, $C_q=\pi(\Ct_q)$ is
irreducible. By \eqref{eq:newton} the form of top degree of $p_q(A,B)$ is
$B^q$, so that of $F_q$ is $A^{2q}$, coming from the term $p_q(B,A)^2$. Hence
$F_q$ has degree $2q$, and $C_q$ has the single point $[0:1:0]$ at infinity.

The vertical maps in \eqref{eq:tower} have degree $6$ and the horizontal ones
degree $q^2$, so $\Ct_q\to C_q$ has degree $6$.

A point of $\Ct_q=\Dt\times_{T,[q]}T$ is a pair $(\Pt,t)$, and $\sigma\in S_3$
fixes it if and only if $\sigma\Pt=\Pt$ and $\sigma t=t$. Its stabilizer is
therefore contained in that of $\Pt$, which by Proposition~\ref{prop:base} has
order at most $2$; so only the twelve points $\Pt$ of ramification of
$\Dt\to\mathbf{P}^1_b$ can carry a ramification point of $\Ct_q\to C_q$, and
there the ramification index is $2$.

Fix such a $\Pt$, with stabilizer $\langle(i\,j)\rangle$, and put
$c=e_i(\Pt)=e_j(\Pt)$; then $e_k(\Pt)=c^{-2}$. The fiber of $\Ct_q\to\Dt$ over
$\Pt$ is a torsor under $T[q]$, of $q^2$ points, and it contains a point with
$t_i=t_j$: take $t_i=t_j=z$ with $z^q=c$ and $t_k=z^{-2}$. Any other point of the
fiber is obtained from it by translation by some $\zeta\in T[q]$, and satisfies
$t_i=t_j$ if and only if $\zeta_i=\zeta_j$, that is
$\zeta_i=\zeta_j=z'$, $\zeta_k=z'^{-2}$ with $z'\in\mu_q$. There are $q$ such
$\zeta$. Hence each of the twelve points carries exactly $q$ ramification points,
and there are $12q$ in all.
\end{proof}

\begin{proof}[Proof of Theorem~\ref{thmA}]
All but the genus is Proposition~\ref{prop:descent}. Riemann--Hurwitz for
$\Ct_q\to C_q$, with Proposition~\ref{prop:kummer} on the left and
Proposition~\ref{prop:descent} on the right, gives
\[
  6q^2-6q=2g(\Ct_q)-2=6\bigl(2g(C_q)-2\bigr)+12q ,
\]
so $2g(C_q)-2=q^2-3q$.
\end{proof}

\begin{remark}\label{rem:toric}
Two natural shortcuts fail. The genus of a curve in a torus is bounded above by
the number of interior lattice points of its Newton polygon, with equality for
nondegenerate curves \cite{Khovanskii1977,Khovanskii1978}; here that bound is
$(q-1)^2$ for $q\not\equiv0\pmod 3$, and $4$ for $\Dt$, whose genus is $1$. So for
$q\ge2$ with $3\nmid q$ neither curve is nondegenerate, and the bound is not
attained. (For $q=3$ the bound is $1=g(C_3)$, and for $q=1$ it is $0=g(C_1)$.) Nor
is $C_4$ the Fermat quartic, which also has genus $3$ and is also a cover of the
line of degree $16$. For the Fermat quartic that cover is Galois, with group
$(\intZ/4)^2$. For $C_4$ the cover $\Phi_4\colon C_4\to D$ is not Galois, since
$\operatorname{Aut}(C_4)$ has order $2$, while the automorphism group of the
Fermat quartic has order $96$ \cite{KuribayashiKomiya1979}. The involution
$\sigma$ of \S\ref{sec:c4} shows that the order is at least $2$. For the upper
bound we computed the automorphism group of the reduction of $C_4$ modulo $5$,
a prime of good reduction, exactly in SageMath; reduction maps
$\operatorname{Aut}(C_4)$ injectively into it.
\end{remark}

\section{Consequences for \texorpdfstring{$n(b)$}{n(b)}}\label{sec:consequences}

\begin{proof}[Proof of Theorem~\ref{thmB}(i)--(iii)]
By Theorem~\ref{thmA}, $g(C_2)=\binom12=0$, $g(C_3)=\binom22=1$ and
$g(C_q)=\binom{q-1}{2}\ge3$ for $q\ge4$.

The assertions in (i) and (ii) about $n(b)$ are proved in \cite{KimuraUmemoto}:
the $b$ with $2\mid n(b)$ are those for which $\eta_b$ is a square, and they are
parametrized by the solutions of $D^2-3E^2=1$; while $3\mid n(b)$ says that
$\eta_b$ is the cube of a unit, which happens only for $b=9$.

Let $q\ge4$. By Theorem~\ref{thmA} the curve $C_q$ is geometrically irreducible
over $\ratQ$ and the genus of its smooth projective model is at least $2$, so
that model has only finitely many rational points \cite{MR718935}; in particular
$C_q$ has only finitely many points with both coordinates in $\intZ$. By
Proposition~\ref{prop:dictionary} every $b$ with $q\mid n(b)$ yields such a point,
and distinct $b$ yield distinct points, since $b$ is recovered from $(A,B)$ by
\eqref{eq:fibre}. Hence there are only finitely many such $b$.
\end{proof}

\begin{remark}
The three regimes match the methods available in each. A Pell equation
parametrizes a conic, and a cubic Thue equation is the search for integral points
on a curve of genus $1$; these are the tools used in \cite{KimuraUmemoto} for
$q=2$ and $q=3$, and Theorem~\ref{thmA} identifies the genus that makes each of
them the right one. For $q\ge4$ no such tool is available, and the conclusion
obtained from Faltings' theorem is not effective.
\end{remark}

Theorem~\ref{thmB}(iii) concerns one $q$ at a time. To pass from it to a statement
about $n(b)$ itself, one must also bound $n(b)$. The bound of
\S\ref{sec:sweep} grows with $|b|$; it suffices for the computation behind
Theorem~\ref{thmC}. A bound independent of $b$ follows from the method of
\S\ref{sec:isopow}. We prove it in \S\ref{sec:isocons}
(Proposition~\ref{iso:uniform}) and deduce Theorem~\ref{thmB}(iv) there.

\section{The sweep}\label{sec:sweep}

Two elementary bounds confine the search. Both concern an arbitrary unit $e>1$ of
a cubic field $K\subset\mathbb{R}$ that is not totally real, with characteristic
polynomial $X^3-BX^2+AX-1$.

\begin{lemma}\label{lem:range}
$|A|\le2e^{1/2}+1$ and $|B-e|\le2e^{-1/2}<2$.
\end{lemma}

\begin{proof}
The conjugates of $e$ are $e$ and a pair $\sigma,\bar\sigma$, and
$\Norm(e)=e|\sigma|^2=1$, so $\sigma=e^{-1/2}\mathrm{e}^{\mathrm{i}\phi}$ for some
$\phi$. Then $B=e+2e^{-1/2}\cos\phi$ and $A=2e^{1/2}\cos\phi+e^{-1}$.
\end{proof}

\begin{lemma}\label{lem:smyth}
$e\ge\rho$, where $\rho=1.3247\ldots$ is the real root of $x^3-x-1$.
\end{lemma}

\begin{proof}
The characteristic polynomial $P=X^3-BX^2+AX-1$ is irreducible, since
$e\neq\pm1$ and $[K:\ratQ]=3$. Its reciprocal is
$X^3P(1/X)=-(X^3-AX^2+BX-1)$, which is $\pm P$ only if $A=B$; but then
$P(1)=A-B=0$, contradicting irreducibility. So $e$ is not reciprocal. Its Mahler
measure is $e$, the other two conjugates having absolute value $e^{-1/2}<1$ by
the proof of Lemma~\ref{lem:range}, so \cite[Theorem]{Smyth1971} applies.
\end{proof}

Applied to $e=\eta_b$, whose characteristic polynomial is \eqref{eq:etamin},
Lemma~\ref{lem:range} gives
\begin{equation}\label{eq:etasize}
  \eta_b<3(b^2+b+1)+2 .
\end{equation}

\begin{corollary}\label{cor:qmax}
Let $|b|\le\beta$. Then $n(b)\le\log(3\beta^2+3\beta+5)/\log\rho$; for
$\beta=10^{30}$, $n(b)\le495$.
\end{corollary}

\begin{proof}
$\eta_b=\eta_0^{\,n(b)}\ge\rho^{\,n(b)}$ by Lemma~\ref{lem:smyth}, while
\eqref{eq:etasize} gives $\eta_b<3\beta^2+3\beta+5$.
\end{proof}

\begin{remark}
The bound of Lemma~\ref{lem:smyth} is attained in the family. For $b=3$ the
fundamental unit of $K_3$ is $\rho$ itself. Indeed $\rho^{13}$ has characteristic
polynomial $X^3-39X^2+12X-1$, which is \eqref{eq:etamin} at $b=3$, and both
$\rho^{13}$ and $\eta_3$ are its real root; so $\eta_3=\rho^{13}$ and
$K_3=\ratQ(\rho^{13})=\ratQ(\rho)$. Thus $\rho$ is a unit $>1$ of $K_3$, and
$\eta_0\ge\rho$ by Lemma~\ref{lem:smyth}; hence $\eta_0=\rho$ and $n(3)=13$.
Read the other way, $n(b)\ge13$ forces
$\rho^{13}\le\eta_b<3(b^2+b+1)+2$ by \eqref{eq:etasize}, and since
$\rho^{13}=38.69\ldots$ this requires $b\ge3$ or $b\le-4$; the value $n(3)=13$
of Theorem~\ref{thmC} sits at that threshold.
\end{remark}

Combining, if $|b|\le\beta$ and $q\mid n(b)$ then the point $(A,B)$ of
Proposition~\ref{prop:dictionary} satisfies
\begin{equation}\label{eq:box}
  |A|\le2\bigl\lfloor(4\beta^2)^{1/2q}\bigr\rfloor+3,
  \qquad -2\le B\le\bigl\lfloor(4\beta^2)^{1/q}\bigr\rfloor+3 .
\end{equation}
Indeed, $\eta_b=u^q$ and $\eta_b<3\beta^2+3\beta+5\le4\beta^2$ for $\beta\ge5$,
so $u<(4\beta^2)^{1/q}$, and Lemma~\ref{lem:range} applied to $e=u$ gives
\eqref{eq:box}. For $\beta=10^{30}$, Corollary~\ref{cor:qmax} gives $q\le495$.
We ran the sweep up to $q=496$, the bound $\log(4\beta^2)/\log\rho$, one
exponent more than needed. The bounds in \eqref{eq:box} are integer parts of
roots of the integer $4\beta^2$, which we compute exactly in integer arithmetic,
so no rounding enters the box.

We obtain Theorem~\ref{thmC} by examining \eqref{eq:box} for each $q$ in turn.
We used two procedures, which agree on the range of $q$ where both are affordable.
The first scans $A$ alone, discarding those $A$ for which $F_q(A,\cdot)$ has no
root modulo one of the forty primes $5\le p\le181$ and factoring the survivors over $\ratQ$; its
cost is governed by the $A$-range, of size $\asymp(4\beta^2)^{1/2q}$. The second
scans the whole box \eqref{eq:box} and evaluates the two power sums at integers
by the recursion \eqref{eq:newton}, so that no polynomial arithmetic occurs and
the dependence on $q$ is only through the length of the recursion. At
$\beta=10^{30}$ we exclude $q=2$ and $q=3$ from the sweep, because their
$A$-ranges have sizes $5.7\cdot10^{15}$ and $5.0\cdot10^{10}$; for them
Theorem~\ref{thmB}(i),(ii) applies instead. At $\beta=10^{12}$ we swept both as
well, and the sweep reproduces those two theorems.

\section{Rational points on \texorpdfstring{$C_4$}{C\_4}}\label{sec:c4}

Theorem~\ref{thmB}(iii) is not effective, and already for $q=4$ it leaves the set
of $b$ with $4\mid n(b)$ undetermined. In this section we determine this set. The main step
is Theorem~\ref{c4:points}, which finds $C_4(\ratQ)$; Theorem~\ref{thmD}
follows in \S\ref{sec:tob}. Several of the statements below rest on exact
computation by machine, and in each case the proof names the computation.
We make precise in \S\ref{rem:certificate} and \S\ref{sec:precision} what is
proved and what is not. In \S\ref{sec:elllog} we determine the integral solutions of
$F_4(A,B)=0$ by a different method, which gives a second proof of Theorem~\ref{thmD}.

For a prime number $\ell$ we write $\mathbf{F}_\ell$ for the field with $\ell$
elements, $\intZ_\ell$ and $\ratQ_\ell$ for the ring of $\ell$-adic integers and the field of
$\ell$-adic numbers,
$\overline{\ratQ}_\ell$ for an algebraic closure of $\ratQ_\ell$, and
$\intZ_{(\ell)}$ for the localization of $\intZ$ at $\ell$.

\subsection{The curve and its quotient}\label{sec:c4curve}
By Theorem~\ref{thmA}, $g(C_4)=3$. The plane model \eqref{eq:Cq} has degree $8$
and is singular; the canonical model of $C_4$ is a smooth plane quartic, which in
suitable coordinates $(t,w)$ is
\begin{equation}\label{eq:quartic}
  w^4-6(t-1)^2w^2+27\,t\,(t^3-4)=0 ,
\end{equation}
of good reduction outside $\{2,3\}$.

\begin{lemma}\label{c4:quotient}
The map $\sigma:w\mapsto-w$ is an involution of $C_4$. Its fixed points are the
four points with $w=0$, and the quotient $C_4/\sigma$ has genus $1$. Explicitly,
$y:=w^2-3(t-1)^2$ satisfies
\begin{equation}\label{eq:H}
  H:\quad y^2=-18t^4-36t^3+54t^2+72t+9 .
\end{equation}
\end{lemma}

\begin{proof}
Only even powers of $w$ occur in \eqref{eq:quartic}, so $\sigma$ is an
involution, and its fixed points are the four points with $w=0$;
Riemann--Hurwitz gives genus $1$ for the quotient. Since
$w^4-6(t-1)^2w^2=y^2-9(t-1)^4$, \eqref{eq:quartic} becomes
$y^2=9(t-1)^4-27t(t^3-4)$, which is \eqref{eq:H}.
\end{proof}

The curve $H$ has the rational points $(t,y)=(0,\pm3)$, and $(0,-3)$ is the
image of $(0,0)$. The map
\begin{equation}\label{eq:HtoE}
  (t,y)\longmapsto\Bigl(\frac{6t^2+24t+6-2y}{3t^2},\
    \frac{36t^4+36t^3+(72-24y)t+2(y-3)^2}{9t^3}\Bigr)
\end{equation}
extends to an isomorphism from $H$ onto
\begin{equation}\label{eq:E}
  E:\quad y^2=x^3-36x,\qquad
  E(\ratQ)=\intZ P_0\oplus(\intZ/2)^2,\quad P_0=(-3,9),
\end{equation}
that sends $(0,-3)$ to the origin. The curve $E$ is the congruent number curve
of $6$, of rank $1$; for $E(\ratQ)$ see \cite[curve 576.c3]{LMFDB}. Write $\varphi\colon C_4\to E$ for the composite of
$C_4\to H$ with \eqref{eq:HtoE}; it is a map of degree $2$.
Put $s:=3(t-1)^2+y$, so that $s=w^2$ on
$C_4$. A point of $H$ lifts to $C_4$ exactly when $s$ is a square in $\ratQ$.

Let $\mathcal{P}$ be the Prym variety of $\varphi$, an abelian surface with
$\Jac(C_4)\sim E\times\mathcal{P}$ over $\ratQ$.

\begin{proposition}\label{c4:rank}
$\operatorname{rank}\mathcal{P}(\ratQ)=2$, and hence
\begin{equation}\label{eq:rank}
  \operatorname{rank}\Jac(C_4)(\ratQ)=3=g(C_4) .
\end{equation}
\end{proposition}

\begin{proof}
For the lower bound, the classes
$D(t,w)=[(t,w)-(t,-w)]$ attached to the points \eqref{eq:c4points} are
anti-invariant under $\sigma$ and so lie in $\mathcal{P}(\ratQ)$. We show that
$D_1:=D(1,3)$ and $D_2:=D\bigl(-\frac1{11},\frac{15}{11}\bigr)$ are independent
modulo torsion. Let $\Lambda\subset\intZ^2$ be the lattice of pairs $(e,f)$ for
which $eD_1+fD_2$ is torsion. This lattice is saturated: if $n(e,f)\in\Lambda$
with $n\neq0$, then $(e,f)\in\Lambda$. So if $\Lambda\neq0$, it contains a
primitive vector. Let $p$ run over the eight primes $5\le p<32$ other than
$11$; they are primes of good reduction, and $11$ is left out because it divides
a denominator of $D_2$. The torsion of $\Jac(C_4)(\ratQ)$ injects into
$\Jac(C_4)(\mathbf F_p)$ for each such $p$, and the orders of these groups have
greatest common divisor $16$, so the torsion has order dividing $16$. Hence
every $(e,f)\in\Lambda$ satisfies $e(16\bar D_1)+f(16\bar D_2)=0$ in
$\Jac(C_4)(\mathbf F_p)$, where the bar denotes reduction modulo $p$. That is,
$\Lambda$ lies in the lattice $L_p$ of these relations, which we computed
exactly in the finite group. The intersection of the eight lattices $L_p$ lies
in $15\,\intZ^2$, so it contains no primitive vector. Therefore $\Lambda=0$.

For the upper bound, write \eqref{eq:quartic} as $w^4-hw^2+fg$ with $f,g,h$
binary quadratic forms; here $h=6(t-1)^2$ and $fg=27t(t^3-4)$, and in that shape
the Prym is the Jacobian of an explicit genus-$2$ curve
\cite[Thm.~1.1]{RitzenthalerRomagny2018}. Over $\ratQ$ no such factorization of
$fg$ exists, because its four roots admit no Galois-stable partition into two
pairs. The smallest field over which such a factorization exists is
$K=\ratQ(2^{1/3})$, which contains the real cube root of $4$. With
$\omega=2^{2/3}$, the curve is
\begin{equation}\label{eq:X}
\begin{split}
  X:\ y^2={}&(4\omega^2+6\omega+8)x^6-(14\omega^2+20\omega+32)x^5\\
            &+(8\omega^2+2\omega+20)x^4+(18\omega^2+28\omega+32)x^3\\
            &-(11\omega^2+8\omega+2)x^2-(4\omega^2+10\omega+28)x-\omega^2-2\omega .
\end{split}
\end{equation}
A $2$-descent by Stoll's algorithm \cite{Stoll2001} gives
$\dim_{\mathbf{F}_2}\mathrm{Sel}^2(\Jac(X)/K)=5$. The sextic factors over $K$
into three irreducible quadratics, so $\dim_{\mathbf F_2}\Jac(X)(K)[2]=2$ and
$\operatorname{rank}\Jac(X)(K)\le5-2=3$, and hence
$\operatorname{rank}\mathcal{P}(K)\le3$. Write
$\mathcal{P}(K)\otimes\ratQ=(\mathcal{P}(\ratQ)\otimes\ratQ)\oplus V$ with $V$
the kernel of the trace, so that
$\operatorname{rank}\mathcal{P}(K)=\operatorname{rank}\mathcal{P}(\ratQ)+\dim V$.
Then $\dim V\ge1$. The point of $C_4(K)$ with $t=(a^2-4a-6)/11$, $a^3=2$, gives a
class $D$ with rational trace, so $3D-\operatorname{Tr}(D)$ lies in $V$. This
class is not torsion: $31$ splits completely in $K$, and its reduction modulo a
prime above $31$ has order $118$, which does not divide $16$.
Here $16$ bounds the torsion of $\Jac(C_4)(K)$, not only of $\mathcal P(\ratQ)$:
each of the good primes $p=5,17,23,29,31$ has a prime of degree $1$ in $K$, the
torsion injects into $\Jac(C_4)(\mathbf F_p)$ by reduction modulo that prime,
and the orders of these groups have greatest common divisor $16$.
Hence $\operatorname{rank}\mathcal{P}(\ratQ)\le3-1=2$.

Since $\Jac(C_4)\sim E\times\mathcal P$ and $E(\ratQ)$ has rank $1$ by
\eqref{eq:E}, \eqref{eq:rank} follows.
\end{proof}

Chabauty's method is therefore unavailable. The quadratic Chabauty method of
Balakrishnan--Dogra \cite{BalakrishnanDogra2018,BalakrishnanDogra2021} applies
instead. It asks for three things \cite[Theorem~1.2(ii)]{BalakrishnanDogra2018}:
equality in \eqref{eq:rank}; a N\'eron--Severi rank of $\Jac(C_4)$ over
$\ratQ$ greater than $1$; and, at the chosen prime $p$, a closure of
$\Jac(C_4)(\ratQ)$ in $\Jac(C_4)(\ratQ_p)$ of finite index (see also
\cite[\S1.4]{BDMTV2019}). The second holds because $\Jac(C_4)$ is isogenous
over $\ratQ$ to $E\times\mathcal P$, and the third is Lemma~\ref{c4:index}.

\subsection{The quadratic Chabauty set}\label{sec:qc}
Let $p=13$, a prime of good reduction for \eqref{eq:quartic}. Take a basis
$\omega_0,\omega_1,\omega_2$ of $H^0(C_{4,\ratQ_p},\Omega^1)$ with $\omega_0$
pulled back from $E$ and $\omega_1,\omega_2$ anti-invariant under $\sigma$.

\begin{lemma}\label{c4:index}
The closure of $\Jac(C_4)(\ratQ)$ in $\Jac(C_4)(\ratQ_p)$ has finite index.
\end{lemma}

\begin{proof}
It suffices that the $p$-adic abelian logarithms of the classes
$\varphi^*(P_0)$, $D(1,3)$ and $D\bigl(-\frac1{11},\frac{15}{11}\bigr)$ span the
dual of $H^0(C_{4,\ratQ_p},\Omega^1)$. The logarithm of $\varphi^*(P_0)$ is
$(c,0,0)$ with $c$ a nonzero multiple of $\log_E(P_0)$, which is nonzero because
$P_0$ has infinite order. The other two classes are anti-invariant, so their
first coordinates vanish, and the $2\times2$ matrix of their remaining
coordinates has determinant of $p$-adic valuation $2$, computed to precision
$p^{9}$.
\end{proof}

The method needs a base point and a correspondence. As base point we take the
rational point $z_0:=(1,3)$. For the correspondence we take $Z$ in
$\ratQ[\sigma]$: up to a constant factor, its endomorphism of
$H^1_{\mathrm{dR}}(C_4)$ is $1+3\sigma^*$. This $Z$ is nice in the sense of
\cite[\S2.3]{BDMTV2019}. In particular its trace is $0$, because $\sigma^*$ acts
as $+1$ on the two-dimensional part pulled back from $E$ and as $-1$ on the
four-dimensional anti-invariant part. With these choices the method attaches to
$C_4$ a function $\Theta\colon C_4(\ratQ_p)\to\ratQ_p$, locally analytic on each
residue disc and built from $p$-adic heights and Coleman integrals, whose value
at a rational point is the sum of the local heights of that point at the primes
other than $p$. Only $\ell=2,3$ contribute, the local height vanishing at primes
of good reduction \cite[Lemma~3.2]{BDMTV2019}.

For $\ell\in\{2,3\}$ write $h_\ell(z)$ for the local height of $z\in C_4(\ratQ_\ell)$.
By Betts and Dogra, as stated in \cite[\S3.1 and Theorem~3.2]{BDMTV2023},
$h_\ell=c_\ell\,j\circ r_\Gamma$. Here $c_\ell$ is a constant, and $r_\Gamma$ is
the retraction of the Berkovich analytification of $C_4$ over $\mathbf{C}_\ell$
onto its minimal skeleton $\Gamma$. The function $j$ is the unique piecewise
polynomial function on $\Gamma$ that vanishes at $r_\Gamma(z_0)$ and has a
prescribed Laplacian $\mu$, built from the graph and the endomorphism induced by
$Z$. In \cite{BDMTV2023} the graph is the dual graph of a regular semistable
model, and the map is reduction to it. We take the stable model with its nodes
resolved; its dual graph is a subdivision of $\Gamma$. A point of
$C_4(\ratQ_\ell)$ reduces to a smooth point $\xi$ of a component. The formal
fiber of $\xi$ is an open ball \cite[Theorem~4.6(2)]{BakerPayneRabinoff2013} and
a connected component of the complement of the vertex set
\cite[Corollary~4.7]{BakerPayneRabinoff2013}, whose closure adds the vertex of
that component. So $r_\Gamma$ maps the point to that vertex, which is its
reduction.

\begin{lemma}\label{c4:hl}
Let $\ell\in\{2,3\}$.
\begin{enumerate}
\item $h_\ell$ is constant on the points of $C_4(\ratQ_\ell)$ that reduce to one
  component of the special fiber of a regular model over $\intZ_\ell$.
\item $h_\ell\circ\sigma=h_\ell$.
\end{enumerate}
\end{lemma}

\begin{proof}
(i) Such a point reduces to a smooth point of a component of multiplicity $1$.
Apply \cite[Theorem~4.6]{BakerPayneRabinoff2013} to the open part of the model
where the special fiber is reduced, base-changed to the ring of integers of
$\mathbf{C}_\ell$. The residue disc of the point is an open ball, and its
closure adds the type-$2$ point $\zeta$ of the component. An open ball meets
no point of the minimal skeleton \cite[Theorem~4.22(1)]{BakerPayneRabinoff2013}.
If $\zeta\in\Gamma$, then $r_\Gamma$ maps every such disc to $\zeta$
\cite[Definition~3.7]{BakerPayneRabinoff2013}. If not, the discs of the
component together with $\zeta$ form a connected set disjoint from $\Gamma$,
on which $r_\Gamma$ is constant.

(ii) The skeleton is intrinsic, so $\sigma$ induces
an isometry of $\Gamma$ that commutes with $r_\Gamma$. Since $Z$ lies in
$\ratQ[\sigma]$, it commutes with $\sigma$, so the Laplacian $\mu$ is
$\sigma$-invariant. Hence $j\circ\sigma-j$ has
Laplacian $0$ and is a constant $k$, and $\sigma^2=1$ gives
$j=j\circ\sigma\circ\sigma=j+2k$, so $k=0$.
\end{proof}

\begin{corollary}\label{c4:upsilon}
For $\ell\in\{2,3\}$, $h_\ell$ takes at most two values on $C_4(\ratQ_\ell)$:
$0$, and the value $v_\ell$ that it takes at the rational points with
$t\in\{-2,0\}$ if $\ell=2$ and with $t\in\{-1,0\}$ if $\ell=3$. Hence
$\Theta(z)\in\Upsilon$ for every $z\in C_4(\ratQ)$, where
\[
  \Upsilon:=\{0,\,v_2,\,v_3,\,v_2+v_3\}\subset\ratQ_p .
\]
\end{corollary}

\begin{proof}
Take the regular models of \S\ref{rem:certificate}. At $\ell=2$ the special
fiber has exactly three components of multiplicity $1$. At $\ell=3$ it has five,
of which exactly three are defined over $\mathbf{F}_3$, and only those can carry
the reduction of a point of $C_4(\ratQ_3)$. In both cases $\sigma$ exchanges two
of the three and fixes the third, and $z_0$ lies on one of the exchanged pair. By
Lemma~\ref{c4:hl}, $h_\ell$ is $0$ on the exchanged pair and takes a single
value $v_\ell$ on the third component. Rational points lie on that
component: those with $t\in\{-2,0\}$ at $\ell=2$, and those with $t\in\{-1,0\}$
at $\ell=3$. So $v_\ell$ is the local height at such a point.
\end{proof}

Put
\[
  \mathcal{S}\ :=\ \{z\in C_4(\ratQ_p)\ :\ \Theta(z)\in\Upsilon\} .
\]
By Corollary~\ref{c4:upsilon}, $\mathcal S$ contains $C_4(\ratQ)$. It is finite.
Indeed, the method yields a quadratic Chabauty pair \cite[Lemma~3.7]{BDMTV2019},
so by property~(a) of \cite[\S1.4]{BDMTV2019} only finitely many points of
$C_4(\ratQ_p)$ give each value of $\Theta$, and $\Upsilon$ is finite.

The four elements of $\Upsilon$ are distinct. From the regular models of
\S\ref{rem:certificate} we obtain $v_2=2106\log_{13}2$ and
$v_3=-4212\log_{13}3$, in the normalization of $Z$ that we use. If two elements
of $\Upsilon$ coincided, then $\log_{13}(2^e3^f)=0$ for some integers
$(e,f)\neq(0,0)$. But the kernel of $\log_{13}$ on $\ratQ_{13}^\times$ lies in
the subgroup generated by $13$ and the roots of unity. The only roots of unity
in $\ratQ$ are $\pm1$, so a rational number in this kernel is of the form
$\pm13^k$. The number $2^e3^f$ is not of this form.

\begin{proposition}\label{c4:S}
The set $\mathcal{S}$ consists of $124$ simple solutions defined over $\ratQ_p$,
the point $(0,0)$, and possibly one further point $z_1\in C_4(\ratQ_p)$ on the
residue disc of $(0,0)$. Of the $124$, $120$ lie on the twenty residue discs not
containing $(0,0)$ and four on the disc of $(0,0)$.
\end{proposition}

\begin{proof}
There are twenty-one residue discs,
and on each of them $\Theta$ takes each value of $\Upsilon$
twice over $\overline{\ratQ}_p$, so that $168$ zeros are counted before any is
tested for rationality. Of these, $124$ are simple zeros defined over $\ratQ_p$:
$120$ on the twenty discs not containing the rational point $(0,0)$, and four on
the disc of $(0,0)$, at two of the four values of $\Upsilon$. At a third value
there is no zero over $\ratQ_p$ on that disc. At the fourth, $v_0=\Theta(0,0)$,
the computation does not separate the two zeros on that disc: one is $(0,0)$,
and the other either coincides with it or is a point $z_1\in C_4(\ratQ_p)$ close
to it.

The search that finds these zeros covers all of $\intZ_p$. It splits a disc
into its $p$ residue classes and discards a class only where the series takes
unit values. It accepts a zero only on a class on which the Newton polygon shows
exactly one zero; by Weierstrass preparation, that zero lies in $\ratQ_p$ and is
simple. Every class is settled in this way except one, the class of $(0,0)$ at
the value $v_0$. No zero is lost beyond the working precision, by Lemma~\ref{c4:tail}.

We carried out the computation in an independent SageMath implementation
following \cite{BDMTV2019}, with Coleman integration on the plane quartic by the
algorithm of Tuitman \cite{Tuitman2016,BalakrishnanTuitman2020}.
\end{proof}

Two internal checks are worth recording, because neither is automatic. First,
the count of zeros is unchanged when the working precision is raised from
$p^{4}$ to $p^{9}$ and again to $p^{13}$, although every intermediate quantity
differs between the runs. Second, $\sigma$ permutes the $124$ simple solutions in
pairs, although for the $120$ off the disc of $(0,0)$ the two discs of a pair
are assembled from different Teichm\"uller points and different Frobenius data.
We test not the parity of $124$ but that the computed solutions match up under
$\sigma$ point by point. On the disc of $(0,0)$, moreover, the computed
expansion of $\Theta-v_0$ is an even function of the uniformizer to the working
precision, as $\sigma$-invariance of $\Theta$ would make it. This evenness would
give $z_1=(0,0)$, but we do not use it.

Eight of the $124$ are rational points \eqref{eq:c4points}
below, and $(0,0)$ is the ninth; $116$ are not accounted for.

\subsection{The sieve}\label{sec:sieve}
The remaining $116$ points are removed by a Mordell--Weil sieve
\cite{BruinStoll2010}, applied to $E$ rather than to $\Jac(C_4)$.
A point $z\in C_4(\ratQ)$ has
$\varphi(z)\in E(\ratQ)=\intZ P_0\oplus(\intZ/2)^2$ by \eqref{eq:E}, so
$\varphi(z)=m(z)P_0+\tau(z)$ with $m(z)\in\intZ$ and
$\tau(z)\in E(\ratQ)_{\mathrm{tors}}$ uniquely determined by $z$; and $s$ takes
a square value there, being $w^2$ on $C_4$.

\begin{proposition}\label{c4:sieve}
None of the $116$ simple solutions in $\mathcal{S}$ that are not among the
points of \eqref{eq:c4points} is rational.
\end{proposition}

\begin{proof}
We sieve the condition above on the pair $(m,\tau)$ modulo the orders of $P_0$
in $E(\mathbf{F}_\ell)$ for a large set of auxiliary primes $\ell$. For an
explicit modulus $M$ of $46$ digits, the sieve leaves
\begin{equation}\label{eq:sieveout}
  m\equiv0,-1\ \ \text{or}\ \ m\equiv0,-1,-2 \pmod M,
  \qquad v_{13}(M)=2 ,
\end{equation}
for $\tau=O$ and for $\tau=(6,0)$ respectively, and no $m$ at all for
$\tau=(0,0)$ and $\tau=(-6,0)$. On the other side, every $z\in\mathcal{S}$ has a
$13$-adic multiplier
\[
  m(z)=\log_{13}\varphi(z)\big/\log_{13}P_0\ \in\ \ratQ_{13} ,
\]
which is well defined because the torsion of $E(\ratQ)$ has order prime to $13$.
For rational $z$ this is the integer above. So $z$ is not rational unless $m(z)\in\intZ_{13}$ and is
admissible modulo $13^{2}$ by \eqref{eq:sieveout}. That excludes $112$ of the
$116$. For the four survivors we raise the $13$-part of $M$ to $13^4$ and run
the sieve again on their residue classes alone; this excludes all four. Primes $\ell$ with
$13^4\mid\operatorname{ord}(P_0\bmod\ell)$ and small cofactor are available
because $E$ has complex multiplication.
\end{proof}

\begin{lemma}\label{c4:disc}
The residue disc of $(0,0)$ in $C_4(\ratQ_{13})$ contains no rational point other
than $(0,0)$.
\end{lemma}

\begin{proof}
The involution $\sigma$ maps that disc to itself, and $(0,0)$ is its only fixed
point in $C_4(\ratQ_{13})$: there is none at infinity, and the affine ones have
$w=0$, so $t=0$ or $t^3=4$, while $4$ is not a cube modulo $13$. Suppose
$z\neq(0,0)$ were a rational point on the disc. Then so is $\sigma(z)$, and it
differs from $z$ and from $(0,0)$. The four simple solutions on the disc are
among the $116$ of Proposition~\ref{c4:sieve}, so neither $z$ nor $\sigma(z)$ is
one of them; hence $\Theta(z)=\Theta(\sigma(z))=v_0$ by
Proposition~\ref{c4:S}, and $(0,0)$, $z$, $\sigma(z)$ would be three distinct
zeros of $\Theta-v_0$ on a disc where it has two.
\end{proof}

\begin{theorem}\label{c4:points}
In the coordinates of \eqref{eq:quartic},
\begin{equation}\label{eq:c4points}
  C_4(\ratQ)=\Bigl\{(1,\pm3),\ (-1,\pm3),\ (-2,\pm6),
    \bigl(-\tfrac1{11},\pm\tfrac{15}{11}\bigr),\ (0,0)\Bigr\} .
\end{equation}
\end{theorem}

\begin{proof}
The nine points lie on \eqref{eq:quartic}. Conversely, $C_4(\ratQ)\subset\mathcal{S}$
by Corollary~\ref{c4:upsilon}. By Proposition~\ref{c4:S}, $\mathcal{S}$ consists of
$124$ simple solutions, the point $(0,0)$ and possibly $z_1$. Eight of the $124$
are the points of \eqref{eq:c4points} other than $(0,0)$, and the other $116$ are
not rational by Proposition~\ref{c4:sieve}. If $z_1$ exists, it lies on the disc
of $(0,0)$ and differs from $(0,0)$, so it is not rational by
Lemma~\ref{c4:disc}.
\end{proof}

As a control, the nine points give $\tau=O$ with $m=0$ for $(0,0)$ and $m=-1$
for $t=-2$, and $\tau=(6,0)$ with $m=0$ for $t=1$, $m=-1$ for $t=-1$ and
$m=-2$ for $t=-\frac1{11}$. This agrees with \eqref{eq:sieveout}.

\subsection{Descent to \texorpdfstring{$b$}{b}}\label{sec:tob}
The canonical map from \eqref{eq:Cq} to \eqref{eq:quartic} is given by three
quintic adjoints of the degree-$8$ model, and the fiber over a point of
\eqref{eq:quartic} is cut out on that model by two of them. The nine points of
\eqref{eq:c4points} give
\begin{center}
\begin{tabular}{lll}
\hline
point of \eqref{eq:quartic} & $(A,B)$ & $b$\\
\hline
$(0,0)$ & $(-2,1)$ & $-3$\\
$(-2,6)$ & $(-4,11)$ & $-75$\\
$(-2,-6)$ & $(0,3)$ & $5$\\
$(-1,-3)$ & $(1,13)$ & $96$\\
$(1,3)$ & $(3,3)$ & $0$\\
$\bigl(-\frac1{11},\frac{15}{11}\bigr)$ &
  $\bigl(-\frac{123}{23},\frac{5931}{529}\bigr)$ &
  $-\frac{22168640}{279841}$\\
$(1,-3)$, $(-1,3)$, $\bigl(-\frac1{11},-\frac{15}{11}\bigr)$ &
  \multicolumn{2}{l}{no point off the base locus}\\
\hline
\end{tabular}
\end{center}
where $b$ is recovered from $(A,B)$ by \eqref{eq:fibre}. We checked each row in
both directions: the fiber over the point, cut out by two of the adjoints, gives
$(A,B)$, and the adjoints map $(A,B)$ back to the point. The base locus of the canonical map meets the degree-$8$
model in the two rational points $(A,B)=(1,1)$ and $(-1,-1)$, and $b=0$ at both.

\begin{proof}[Proof of Theorem~\ref{thmD}]
Suppose $4\mid n(b)$. By Proposition~\ref{prop:dictionary} there is
$(A,B)\in\intZ^2$ on $C_4$ satisfying \eqref{eq:fibre}. Either $(A,B)$ lies in
the base locus of the canonical map, in which case $b=0$, which is excluded; or
the adjoints evaluated at $(A,B)$ give a rational point of \eqref{eq:quartic},
which by Theorem~\ref{c4:points} is one of nine, and $b$ is then one of the
values in the table. Of these, $b=-\frac{22168640}{279841}$ is not an integer and
$b=0$ is excluded, leaving $b\in\{-75,-3,5,96\}$. Conversely
$n(-75)=n(5)=n(96)=4$ and $n(-3)=8$.
\end{proof}

The remaining three rational points of $C_4$ lie over the base locus, where
$b=0$: the canonical map separates the branches of the ordinary triple point
$(A,B)=(1,1)$, whose only rational branch gives $(-1,3)$, and of the node
$(-1,-1)$, whose two branches give $(1,-3)$ and
$\bigl(-\frac1{11},-\frac{15}{11}\bigr)$. The involution $\sigma$ does not
preserve $b$: $(-2,6)$ and $(-2,-6)$ give $b=-75$ and $b=5$, and $(-1,-3)$
gives $b=96$ while $(-1,3)$ gives the excluded value $b=0$.

\subsection{What the argument rests on}\label{rem:certificate}
The computations behind Theorem~\ref{c4:points} are of three kinds, and they
have different standing.

First, some are exact. The sieve of Proposition~\ref{c4:sieve} uses integer and
finite-field arithmetic throughout, and a residue class is discarded only when it
has been shown to be empty. Proposition~\ref{c4:rank} is exact in the same sense:
the $2$-descent over $K$, the reduction certificate and the point counts are all
exact, and none of them requires a parity conjecture or the finiteness of a
Tate--Shafarevich group.

Second, the computation of Lemma~\ref{c4:index} and Proposition~\ref{c4:S} is
$p$-adic. It is carried out modulo an explicit power of $p$, and every stage
carries a proved precision. The input on which it turns is a $p$-adic height
pairing of two divisor classes of degree zero. In the series whose zeros are
being located, this pairing is multiplied by a series with coefficients of
valuation $-3$, so three digits of the pairing are consumed before any zero can
be placed. We proved the pairing to $p^{16}$. We recomputed $\mathcal{S}$ with a
certificate of $p^{13}$ and reproduced the count, the table of counts across
$\Upsilon$ and the $\sigma$-symmetry. We located the $116$ remaining points of
Proposition~\ref{c4:sieve} to $p^{10}$, six digits beyond what its sieve
consumes.

Third, and weaker, is the input to Corollary~\ref{c4:upsilon}: the local
heights away from $p$ are computed from regular models of \eqref{eq:quartic}
over $\intZ_2$ and $\intZ_3$, both constructed here. The one
at $3$ comes from the Newton polytope of \eqref{eq:quartic}
\cite{Dokchitser2021}. Its component group, of order $16$, and the reduction of
the points \eqref{eq:c4points} agree with those of the regular model returned by
Magma \cite{Magma}.

At $2$ no chart satisfies the hypothesis of the Newton-polytope construction,
the reduction being a quadruple line. Instead, we implement Liu's
desingularization of double covers of regular surfaces \cite{Liu2026} and run it
twice, on models over $\intZ_2$:
\[
  \mathbf{P}^1_{\intZ_2}\longleftarrow\mathcal{E}\longleftarrow\mathcal{C}_4 .
\]
The first run resolves the Weierstrass model of $E$ over $\mathbf{P}^1_t$ and
gives $\mathcal{E}$. The second resolves the double cover $w^2=s$ of what the
first returns, with $s$ as in \S\ref{sec:c4curve}, and gives the regular model $\mathcal{C}_4$ of \eqref{eq:quartic}. Four
blow-ups in the first round and none in the second give eight components, of
multiplicities $(1,1,1,2,2,2,2,4)$, and component group $\intZ/8$. Contracting
two $(-1)$-curves gives a model with six components whose intersection matrix
agrees, up to relabeling, with that of the regular model returned by Magma
\cite{Magma}. This agreement is a check, not an input to the proof.

Four checks are independent of that agreement. Computed with the geometric
genera of the components, the genus formula gives $3=g(C_4)$; this holds only
if every component is regular. The same code, run on Weierstrass models of
elliptic curves, returns the Kodaira type and the component group that Tate's
algorithm gives, at $2$ and at $3$. Run at $\ell=3$ it returns the component
group, the partition of \eqref{eq:c4points} and the local pairing of the
Newton-polytope model above. And it reproduces Liu's own Example 3.10, a family of curves with integer
parameters $m$ and $n$, statement by statement
for $m$ even, down to the equation $w^2+w+x^3+\bar m=0$ of the elliptic
component. For $m$ odd it does not. In Liu's notation, from equation (3.10)
in paragraph (IV.3.2) of that example we compute $\lambda_{p_3}=3$ where the paper states $4$, and the
resolution then has fourteen components and no elliptic one. Liu has kindly
confirmed our computation for $m$ odd (private communication, November 2026).
The resolution of $C_4$ does not
reach that case. The weak point that remains is the source: \cite{Liu2026} is
not yet published; it is currently under revision.

These four checks concern the fiber. A fifth check concerns the points on the
fiber, from which we read the local heights. The placement of
\eqref{eq:c4points} on the special fiber at $2$ is forced by the intersection
matrix alone. The partition it induces agrees with the one given by the
component of $\varphi(z)-\varphi(z_0)$ in the N\'eron model of $E$ at $2$, for
$z$ in \eqref{eq:c4points}, a computation that uses no model of $C_4$.

\subsection{Precision, and the completeness of
\texorpdfstring{$\mathcal{S}$}{S}}\label{sec:precision}
The remaining question is whether a finite computation can miss a zero.
Balakrishnan, Dogra, M\"uller, Tuitman and Vonk state the corresponding result
for the split Cartan modular curve of level $13$ without reservation, and
rightly so. An elementary estimate \cite[\S6.4.1]{BDMTV2019} shows that the
coefficients $c_n$ of the series in question satisfy
$v_p(c_n)\ge C-2\lfloor\log_pn\rfloor$. Here $C$ is the least valuation of the
scalars entering the series, and the term $2\lfloor\log_pn\rfloor$ is what two
rounds of integration cost. Since $\log_pn=o(n)$, every slope of the Newton
polygon beyond the computed precision then exceeds $-1$. Only a slope at most
$-1$ can belong to a zero in the disc under consideration. So no zero can hide
beyond the computed precision, and the finite computation settles the
question.

\begin{lemma}\label{c4:tail}
On each of the twenty-one residue discs of $C_4(\ratQ_{13})$ and for each
element $v$ of $\Upsilon$, no zero of $\Theta-v$ on the disc is carried by the
coefficients of its expansion beyond the computed range. In particular the
computation of Proposition~\ref{c4:S} finds all of $\mathcal{S}$.
\end{lemma}

\begin{proof}
We carry out the estimate of \cite[\S6.4.1]{BDMTV2019} for
\eqref{eq:quartic} at $p=13$. Three points require attention.

First, the estimate applies to a linear combination of finitely many power
series $F_i$ with coefficients in $\intZ_p$, of their single integrals $\int
F_i$ and of their double integrals $\int F_i(\int F_j)$. The function $\Theta$
is not given in that form, but as a quadratic expression in the
coordinates of the Frobenius trivialization. Expanding it yields a quadratic
form over $\ratQp$ in the six single Coleman integrals, together with one
double integral and one power series. The required form is then recovered from
$(\int F_i)(\int F_j)=\int F_i(\int F_j)+\int F_j(\int F_i)$.

Second, the estimate assumes that the $F_i$ have coefficients in $\intZ_p$.
This needs proof only on the residue disc of $(0,0)$, the only disc on which the
auxiliary polynomial of Tuitman's algorithm vanishes. The local parameter there
is $w$.
Substituting $t=w^2X$, with $X$ a new indeterminate, in \eqref{eq:quartic} and
dividing by $w^2$ gives a polynomial in $X$ and $w$ over $\intZ_{(13)}$ whose
value at $w=0$ is $-6-108X$. This has the root $X=-1/18$, a $13$-adic unit, at
which its derivative in $X$ equals $-108$, again a unit. By Hensel's lemma over
$\intZ_{13}[[w]]$ one has $t=w^2X(w)$ with $X(w)$ a unit of $\intZ_{13}[[w]]$,
and the integrality follows.

Third, on each residue disc the series is written in a local parameter $\xi$
that maps the disc onto $p\intZ_p$. On the twenty discs not containing $(0,0)$,
$\xi=t-t_0$ with $t_0\in\intZ_p$ a lift of the residue of $t$. On the disc of
$(0,0)$, $\xi=w$, and by the second point every $w\in p\intZ_p$ gives a point of
that disc. So the series is read in $s=\xi/p$, whose coefficients $d_n=c_np^n$ satisfy $v_p(d_n)\ge C+n-2\lfloor\log_pn\rfloor$.
This bound is not monotone in $n$, so it is not
immediate that the coefficients beyond the computed range may be disregarded.
For $p=13$ one has $\lfloor\log_pn\rfloor\le(n-1)/4$, which gives the monotone
minorant $C+(n+1)/2$ and makes the verification finite. We check something
stronger than that the tail misses the content of the series: every coefficient
beyond the computed range is negligible at the precision to which the series is
known. This is needed because of how the zeros are located. The search fixes
one $13$-adic digit at a time: it shifts the variable by that digit, rescales it
by $13$, and reads the Newton polygon of the resulting series. No such change
of variable lowers the valuation of a coefficient, so a tail that is negligible
at the working precision remains negligible at every step.
The threshold is therefore the content of the series plus the precision to which
it is known. The content is read off a coefficient whose valuation lies below its
own precision. Such a valuation is exact, and so bounds the true content from
above, as the comparison requires. In the run with certificate $p^{9}$ the
numbers are the same on each of the twenty-one discs and for each of the four
elements of $\Upsilon$: $C=0$, the content is $2$, and the threshold is
$2+9=11$. Since $n-2\lfloor\log_{13}n\rfloor\ge11$ for every $n\ge11$, only
$d_0,\dots,d_{10}$ can matter, and we computed sixty-seven coefficients. The
run with certificate $p^{13}$, which locates the points used by the sieve of
Proposition~\ref{c4:sieve}, gives again the same numbers on every disc and for
every element of $\Upsilon$: $C=0$, the content is $2$, and the threshold is
$2+13=15$. Since $n-2\lfloor\log_{13}n\rfloor\ge15$ for every $n\ge17$, only
$d_0,\dots,d_{16}$ can matter, and we computed eighty-three coefficients. No
zero is therefore carried by the tail in either run.
\end{proof}

One normalization is easily overlooked. In \S\ref{sec:qc} the global $p$-adic
height enters as a symmetric bilinear form $h$ on
$H^0(C_{4,\ratQp},\Omega^1)^*$, written in the coordinates dual to the basis
$\omega_0,\omega_1,\omega_2$ of $H^0(C_{4,\ratQp},\Omega^1)$. Let $E'$ be the
Weierstrass model of $E$ through which $\varphi$ is computed, $\omega_{E'}$ its
invariant differential, $P$ a generator of $E'(\ratQ)$, $\hat h_{E'}$ the
$p$-adic height on $E'$ and $\log_{E'}$ its formal logarithm. The basis is
chosen so that $\omega_0=\Lambda\,\varphi^*\omega_{E'}$ for an explicit
$\Lambda\in\ratQ$, while $\omega_1$ and $\omega_2$ are anti-invariant under
$\sigma$. Put $\kappa=\varphi^*(P)$.

\begin{lemma}\label{c4:norm}
$h(\kappa,\kappa)\big/\bigl(\int_\kappa\omega_0\bigr)^2
=\hat h_{E'}(P)\big/\bigl(2(\Lambda\log_{E'}P)^2\bigr)$.
\end{lemma}

\begin{proof}
The class $\kappa$ is $\sigma$-invariant, so its second and third coordinates
vanish, and the projection formula for $\varphi$ gives
\[
  \int_\kappa\omega_0=\Lambda\deg(\varphi)\log_{E'}(P),\qquad
  h(\kappa,\kappa)=\deg(\varphi)\,\hat h_{E'}(P).
\]
The entry $h(\kappa,\kappa)\big/\bigl(\int_\kappa\omega_0\bigr)^2$ of $h$ is
therefore $\hat h_{E'}(P)/(\Lambda\log_{E'}P)^2$ divided by $\deg\varphi=2$.
\end{proof}

We use Lemma~\ref{c4:norm} in the computation of Proposition~\ref{c4:S}.
The second and third coordinates of $\kappa$ vanish, so the entry of $h$ in the
first coordinate is $h(\kappa,\kappa)\big/\bigl(\int_\kappa\omega_0\bigr)^2$,
and the lemma obtains it from the $p$-adic height on $E'$.

We confirmed the factor $\frac12$ numerically on the bielliptic curve
$y^2=x^6-2x^4-2x^2+4$ of genus $2$. Its height pairing is known independently,
because both quotients are elliptic curves of rank $1$. At nine of its
rational points and at $p=11$ and $p=19$, the height computed by our
implementation equals half the pairing obtained from the two elliptic curves, to
the working precision.

This completes the proof of Theorem~\ref{thmD}. The proof relies at several
points on exact computation carried out by machine: the $2$-descent behind
Proposition~\ref{c4:rank}, the quadratic Chabauty computation of
Lemma~\ref{c4:index} and Proposition~\ref{c4:S}, the
Mordell--Weil sieve of Proposition~\ref{c4:sieve}, and the two regular models of
\S\ref{rem:certificate}.

\subsection{Integral points by elliptic logarithms}\label{sec:elllog}
Theorem~\ref{thmD} needs only the integral solutions of $F_4(A,B)=0$, not all
of $C_4(\ratQ)$. These can also be found by elliptic logarithms, in the manner
of \cite{StroekerTzanakis1994,GebelPethoZimmer1994}, without quadratic
Chabauty.

\begin{remark}\label{rem:elllog}
The integral solutions of $F_4(A,B)=0$ are
\[
  (A,B)=(-2,1),\ (-4,11),\ (0,3),\ (1,13),\ (3,3),\ (1,1),\ (-1,-1).
\]
The first five are those of the table in \S\ref{sec:tob}, and the last two are
the base points of the canonical map, where $b=0$. So this proves
Theorem~\ref{thmD} again, without Theorem~\ref{c4:points} and without the
regular models of \S\ref{rem:certificate}. It does not determine $C_4(\ratQ)$.
We sketch the argument. Below, $x$ and $y$ are the coordinates of $E$ in
\eqref{eq:E}, $\varphi(A,B)$ is the image under $\varphi$ of the point of $C_4$
given by $(A,B)$, and $\hat h$ is the canonical height normalized so that
$\hat h\approx\frac12h(x)$, as in \cite{Silverman1990}.

The smooth projective model of $C_4$ has four points over $A=\infty$. They are
conjugate over $\ratQ$ and defined over $\mathbb{R}$; at them, $B/A^2$ tends to
the four roots of $v^4-16v^3+20v^2-8v+1$. The map $\varphi$ sends each of
them to a point $Q_i$ of $E$ with $2Q_i=R$, where $R=(12,-36)=-P_0+(0,0)$. For
$|A|\ge320$, an integral solution lies on one of the four branches through these
points, by Rouch\'e's theorem. By Schwarz's lemma, the integral $\delta$ of
$dx/(2y)$ along the image of the real arc of that branch, from $Q_i$ to
$\varphi(A,B)$, satisfies $|\delta|\le548/|A|$. Write $\varphi(A,B)=mP_0+\tau$
with $\tau$ a torsion point, as in \S\ref{sec:sieve}, and put
$P_1=P_0+(-6,0)=(18,-72)$, so that $\varphi(A,B)=mP_1+\tau'$ with $\tau'$
torsion and $R=-P_1+(6,0)$. Let $\varpi$ be the real period of $dx/(2y)$, and
$u_1\in(0,\varpi/2)$ the elliptic logarithm of $P_1$; the elliptic logarithm
of $(6,0)$ is $\varpi/2$. Doubling gives
\[
  2\delta=(2m+1)\,u_1+k\,\varpi/2\qquad(k\in\intZ,\ |k|\le|2m+1|+8),
\]
a linear form in two elliptic logarithms of rational points. It does not vanish,
since $P_1$ has infinite order. On the other side, $x\circ\varphi$ is a quotient
of two polynomials in $A,B$ of weighted degree $8$, with weights $1$ and $2$;
on the four branches the denominator does not vanish. For $|A|\ge320$,
Silverman's bound \cite[Example~2.2]{Silverman1990} for the difference between
the Weil height and the canonical height then gives
$m^2\hat h(P_0)\le4\log|A|+9.02$, where $\hat h(P_0)>0.4442$. David's lower
bound \cite[Th\'eor\`eme~2.1]{David1995}, applied with two logarithms on the
same curve, gives $|m|<2^{81}$. An elementary reduction, with one rational
approximation of $2u_1/\varpi$, then gives $|m|\le25$. Finally, we find the
integral solutions with $|A|<320$ by factoring $F_4(A,B)$ for each $A$, and
those with $|m|\le25$ by pulling the points $mP_0+\tau$ back to $C_4$. Both
searches give the seven points above.

Every inequality in this argument is checked in exact rational arithmetic, by
the script \texttt{check-c4-elllog.gp} of \cite{KimuraCode2}. The numbers $u_1$
and $\varpi/2$ are the integrals of $dx/\bigl(2\sqrt{x^3-36x}\bigr)$ from $x=18$
and from $x=6$ to $\infty$. We enclose them between rational numbers by the
duplication formula for such integrals; its proof reduces to a polynomial
identity, which the script checks. The argument rests on the two cited
theorems, on $E(\ratQ)$ as given in \eqref{eq:E}, and on this computation.
\end{remark}

\section{Isomorphisms between the fields \texorpdfstring{$K_b$}{K\_b}}
\label{sec:iso}

Distinct parameters can give the same field: $K_{-1}=K_5$
(Example~\ref{iso:ex}). In this section we ask how far $b\mapsto K_b$ is from
being injective. The question leads to the curves $C_m$ once more, now intersected with
$D$ itself, the plane curve defined by~\eqref{eq:D}. For large $m$ the intersection is controlled by a lower bound for
linear forms in two logarithms, and we handle the exponents left over by an
exact computation.

Field isomorphism problems of this kind have been studied for other
families of cubic fields with one parameter, through cubic Thue equations. For
Shanks's simplest cubic fields $L_m$, the splitting fields of
$X^3-mX^2-(m+3)X-1$, one has $L_m=L_n$ with $-1\le m<n$ only for $m,n$ in an
explicit set of eleven values; this is a theorem of Okazaki, proved by Hoshi
\cite[Theorems~1.1, 1.4 and~1.5]{Hoshi2011}, who matches the isomorphism
classes with the non-trivial solutions of a family of cubic Thue equations.
Hoshi and Miyake \cite{HoshiMiyake2010} give such a correspondence for the
integer specializations $s=m$ of $X^3+sX+s$, where coincidences are frequent:
there are seventeen pairs $(m,n)$ with $-10\le m\le5$ \cite[Corollary~5]{HoshiMiyake2010}. Our
route is different. Since $K_b$ has unit rank one, $K_c=K_b$ makes $\eta_c$ a
power of the fundamental unit of $K_b$ (Lemma~\ref{iso:reform}), and the
question becomes one about powers of a single unit. The Thue-equation methods
decide the question one parameter at a time, and the complete answer for the
simplest cubic fields rests on Okazaki's gap principle for cubic forms of
positive discriminant \cite[\S3]{Hoshi2011}. Our fields are not totally real,
so their cubic forms have negative discriminant. The unit route gives instead
statements uniform in $b$ (Theorem~\ref{iso:main}).

\subsection{Reformulation}\label{sec:isoreform}
We say that a unit $u>1$ of $K_b$ is \emph{of family type} if its characteristic
polynomial $X^3-BX^2+AX-1$ satisfies $(A,B)\in D$, that is, $3B=A^2-3A+9$
by~\eqref{eq:D}. By the definition of $D$, $\eta_b$ is of family type.

\begin{lemma}\label{iso:type}
A unit $u>1$ of $K_b$ is of family type if and only if its characteristic
polynomial is \eqref{eq:etamin} with $b$ replaced by an integer $c$. Then
$c=A/3-1$, $c\neq0,1$, and $u=\eta_c$ as real numbers.
\end{lemma}

\begin{proof}
If $3B=A^2-3A+9$ then $3\mid A^2$, so $3\mid A$. Put $c=A/3-1$. Then $A=3(c+1)$
and $B=(A^2-3A+9)/3=3(c^2+c+1)$. The converse is the definition of $D$. For $c=0$
the polynomial is $(X-1)^3$, while $u\neq1$. For $c=1$ it is
$X^3-9X^2+6X-1$, whose values at $0,\tfrac3{10},1,9$ are
$-1,\tfrac{17}{1000},-3,53$. So it has three real roots, whereas the
characteristic polynomial of $u$ has two non-real roots (\S\ref{sec:setup}).
Finally, for $c\neq0,1$ the polynomial \eqref{eq:etamin} has exactly one real root,
$\eta_c$ (Lemma~\ref{lem:eta}), and $u$ is a real root of it.
\end{proof}

\begin{lemma}\label{iso:reform}
Let $b\neq b'$ be integers different from $0$ and $1$. The following are
equivalent:
\begin{enumerate}
\item $K_b\cong K_{b'}$;
\item $K_b=K_{b'}$ as subfields of $\mathbb{R}$;
\item $\eta_{b'}\in K_b$.
\end{enumerate}
If they hold, let $\eta_0$ be the fundamental unit of $K_b=K_{b'}$. Then
$\eta_b=\eta_0^{\,n(b)}$ and $\eta_{b'}=\eta_0^{\,n(b')}$ with $n(b)\neq n(b')$.
\end{lemma}

\begin{proof}
(i)$\Rightarrow$(ii): Let $\sigma\colon K_{b'}\to K_b$ be an isomorphism.
Composed with the inclusion $K_b\subset\mathbb{R}$, it is a real embedding of
$K_{b'}$, hence the inclusion $K_{b'}\subset\mathbb{R}$, since $K_{b'}$ has only
one. So $\sigma(x)=x$ for every $x\in K_{b'}$.

(ii)$\Rightarrow$(iii) is clear.

(iii)$\Rightarrow$(i): The real root of $f_{b'}$ is
$(b'+1)-\eta_{b'}^{-1}$, which lies in $K_b$. So $K_{b'}\subseteq K_b$, with
equality since both have degree $3$.

The last assertion follows from Lemma~\ref{lem:eta}, applied in the common field.
Moreover $n(b)\neq n(b')$, because $\eta_b\neq\eta_{b'}$: their minimal polynomials
differ in the coefficient $3(b+1)$ of $X$.
\end{proof}

\begin{proposition}\label{iso:bij}
Let $b\neq0,1$ and let $\eta_0$ be the fundamental unit of $K_b$. Then
$c\mapsto n(c)$ is a bijection
\[
  \{c\in\intZ\setminus\{0,1\}:K_c=K_b\}\ \longrightarrow\
  \{j\ge1:\eta_0^{\,j}\text{ is of family type}\}.
\]
\end{proposition}

\begin{proof}
If $K_c=K_b$, then $\eta_c=\eta_0^{\,n(c)}$ is of family type. The map is
injective by the last assertion of Lemma~\ref{iso:reform}. Conversely, suppose
$\eta_0^{\,j}$ is of family type. Lemma~\ref{iso:type} gives $c\neq0,1$ with
$\eta_0^{\,j}=\eta_c$. Then $\eta_c\in K_b$, so $K_c=K_b$ by
Lemma~\ref{iso:reform}, and $n(c)=j$.
\end{proof}

So the fields that coincide with a given $K_b$ correspond to the powers of $\eta_0$
of family type. For $j=n(b)\,m$, the question is whether $\eta_b^{\,m}$ is of
family type. Since $\Phi_m(A,B)$ consists of the coefficients of the $m$-th power
(\S\ref{sec:torus}), this holds if and only if
\begin{equation}\label{eq:Pm}
  P_m(b):=F_m\bigl(3(b+1),3(b^2+b+1)\bigr)=0 ,
\end{equation}
that is, the point of $D$ over $b$ lies on $C_m$.

\begin{example}\label{iso:ex}
The minimal polynomial of $\eta_{-1}$ is $X^3-3X^2-1$, so $(A,B)=(0,3)$. By
\eqref{eq:newton}, $\Phi_2(A,B)=(A^2-2B,\,B^2-2A)$, and $\Phi_4=\Phi_2\circ\Phi_2$.
Hence $\Phi_2(0,3)=(-6,9)$ and $\Phi_4(0,3)=\Phi_2(-6,9)=(18,93)$. So
$\eta_{-1}^{\,4}$ has characteristic polynomial $X^3-93X^2+18X-1$, which is
\eqref{eq:etamin} for
$b=5$. By Lemma~\ref{iso:type}, $\eta_{-1}^{\,4}=\eta_5$, so $K_{-1}=K_5$.
Since $3\nmid-1$ and $(-1)^3-4=-5$ is squarefree, $\eta_{-1}$ is the fundamental
unit by \cite{KimuraUmemoto}. So $n(-1)=1$ and $n(5)=4$, in agreement with
Theorem~\ref{thmC}.
\end{example}

\subsection{Powers of \texorpdfstring{$\eta_b$}{eta\_b}}\label{sec:isopow}

\begin{theorem}\label{iso:main}
Let $b\neq0,1$ and $m\ge2$. Then $\eta_b^{\,m}$ is of family type, that is
$P_m(b)=0$, only for $(b,m)=(-1,4)$.
\end{theorem}

We settle the exponents $2\le m\le30$ by computing $P_m$.

\begin{lemma}\label{iso:small}
For $2\le m\le30$, the only integer root of $P_m$ outside $\{0,1\}$ is $b=-1$
for $m=4$.
\end{lemma}

\begin{proof}
We compute $P_m\in\intZ[b]$ from \eqref{eq:newton} and find its integer roots
exactly. As a check, we also compute $P_m$ as the characteristic polynomial
of the $m$-th power of a root of \eqref{eq:etamin} over $\ratQ(b)$.
\end{proof}

For $m\ge31$ we use a lower bound for linear forms in two logarithms. Fix $b$ and
write $\lambda=\eta_b$ and $B_b=3(b^2+b+1)$. Let $\mu,\bar\mu$
be the non-real conjugates of $\lambda$. Then $|\mu|^2=\lambda^{-1}$, and we may
write $\mu=\lambda^{-1/2}\mathrm{e}^{\mathrm{i}\psi}$ with $0<\psi<\pi$. Put
\[
  \gamma=\mu/\bar\mu=\mathrm{e}^{2\mathrm{i}\psi},\qquad
  \vartheta=\arccos\operatorname{Re}\gamma\in[0,\pi],\qquad
  \mathcal{A}=\{\pm\tfrac\pi3\}+2\pi\intZ ,
\]
so that $2\psi\equiv\pm\vartheta\pmod{2\pi}$. As $B_b=\lambda+2\operatorname{Re}\mu$,
\begin{equation}\label{eq:cosomega}
  \cos\vartheta=\operatorname{Re}\frac{\mu^2}{|\mu|^2}
  =\lambda\bigl(2(\operatorname{Re}\mu)^2-\lambda^{-1}\bigr)
  =\tfrac12\lambda(B_b-\lambda)^2-1 .
\end{equation}

\begin{lemma}\label{iso:identity}
Let $m\ge1$. Then $\eta_b^{\,m}$ is of family type if and only if
$1-2\cos2m\psi=R_m$, where
\[
  R_m=-6\lambda^{-m/2}\cos m\psi+9\lambda^{-m}-2\lambda^{-3m/2}\cos m\psi
      -3\lambda^{-2m}+\lambda^{-3m} .
\]
If $\lambda^m\ge16$, then $|R_m|\le6\lambda^{-m/2}+10\lambda^{-m}$.
\end{lemma}

\begin{proof}
The conjugates of $u=\eta_b^{\,m}$ are $\lambda^m$ and
$\lambda^{-m/2}\mathrm{e}^{\pm\mathrm{i}m\psi}$, and their product is $1$. Hence
$B=\lambda^m+2\lambda^{-m/2}\cos m\psi$ and
$A=\sum_i e_i^{-1}=\lambda^{-m}+2\lambda^{m/2}\cos m\psi$. Expanding, and using
$3-4\cos^2m\psi=1-2\cos2m\psi$,
\[
  3B-(A^2-3A+9)=\lambda^m\bigl((1-2\cos2m\psi)-R_m\bigr).
\]
For the bound, note that $2\lambda^{-3m/2}+3\lambda^{-2m}+\lambda^{-3m}\le
(\tfrac24+\tfrac3{16}+\tfrac1{256})\lambda^{-m}<\lambda^{-m}$ when
$\lambda^m\ge16$.
\end{proof}

\begin{lemma}\label{iso:angle}
Let $x\in\mathbb{R}$ and $\varepsilon=\operatorname{dist}(x,\mathcal{A})$. Then
\[
  \varepsilon/\pi\ \le\ |\cos x-\tfrac12|\ \le\ \varepsilon ,
\]
and if $\varepsilon\le\frac1{10}$, then $|\cos x-\frac12|\ge0.815\,\varepsilon$.
\end{lemma}

\begin{proof}
Let $x'\in[0,\pi]$ with $x'\equiv\pm x\pmod{2\pi}$. Then $\cos x'=\cos x$, and
$\varepsilon=|x'-\frac\pi3|$, the nearest point of $\mathcal{A}$ to $x'$ being
$\frac\pi3$. Now
\[
  \cos x-\tfrac12=\cos x'-\cos\tfrac\pi3
  =-2\sin\tfrac{x'+\pi/3}{2}\,\sin\tfrac{x'-\pi/3}{2} .
\]
Here $|\sin\frac{x'-\pi/3}2|=\sin\frac\varepsilon2$, as
$\frac\varepsilon2\in[0,\frac\pi3]$, and
$\frac{x'+\pi/3}2\in[\frac\pi6,\frac{2\pi}3]$, where $\sin\ge\frac12$. So
$\sin\frac\varepsilon2\le|\cos x-\frac12|\le2\sin\frac\varepsilon2$, and
$\frac2\pi t\le\sin t\le t$ for $0\le t\le\frac\pi2$ gives the first assertion.
If $\varepsilon\le\frac1{10}$, then $\frac{x'+\pi/3}2$ lies within $\frac1{20}$ of
$\frac\pi3$, so $\sin\frac{x'+\pi/3}2\ge\frac{\sqrt3}2-\frac1{20}>0.816$. Also
$\sin\frac\varepsilon2\ge\frac\varepsilon2(1-\frac{\varepsilon^2}{24})
\ge\frac\varepsilon2(1-\frac1{2400})$. Hence
$|\cos x-\frac12|\ge0.816(1-\frac1{2400})\,\varepsilon\ge0.815\,\varepsilon$.
\end{proof}

\begin{lemma}\label{iso:beta}
The number $\gamma$ is not a root of unity. Its degree is even and at most $6$,
and its absolute logarithmic height is $h(\gamma)=\frac12\log\lambda$. In particular
$0<\vartheta<\pi$.
\end{lemma}

\begin{proof}
If $\gamma^n=1$, then $\mu^n=\bar\mu^n$. So two conjugates of $\eta_b^{\,n}$ would
coincide, although $\eta_b^{\,n}\notin\ratQ$ generates the cubic field $K_b$.
Let $L$ be the Galois closure of $K_b$. Since $\disc{K_b}<0$ is not a square,
$\operatorname{Gal}(L/\ratQ)\cong S_3$, acting on the roots $\lambda,\mu,\bar\mu$.
So $\gamma\in L$ has degree at most $6$. The degree is even: the conjugates of
$\gamma$ include $\bar\gamma=\gamma^{-1}$, so they are permuted by $x\mapsto x^{-1}$,
an involution whose only fixed points $\pm1$ are not conjugates of $\gamma$. As
$\sigma$ runs over
$\operatorname{Gal}(L/\ratQ)$, the numbers $\sigma(\gamma)=\sigma(\mu)/\sigma(\bar\mu)$
run over the six quotients $e_i/e_j$, $i\neq j$, of the roots, each once. Only
$\lambda/\mu$ and $\lambda/\bar\mu$ have absolute value greater than $1$, namely
$\lambda^{3/2}$. Since $\gamma$ is a unit,
$h(\gamma)=\frac16\sum_\sigma\log\max(1,|\sigma(\gamma)|)=\frac16\cdot2\cdot\frac32\log\lambda$.
Finally, $\vartheta\in\{0,\pi\}$ would give $\gamma=\pm1$.
\end{proof}

\begin{lemma}\label{iso:linear}
Let $b\neq0,1$ and $m\ge31$, and suppose that $\eta_b^{\,m}$ is of family type.
Then $\lambda>9.3$ if $b\neq-1$, and
$\operatorname{dist}(2m\psi,\mathcal{A})\le4.3\lambda^{-m/2}$. Moreover there is an
odd integer $j$ with $1\le j\le3m$ and
\begin{equation}\label{eq:linform}
  |j\pi-3m\vartheta|\le13\lambda^{-m/2} .
\end{equation}
\end{lemma}

\begin{proof}
By Lemma~\ref{lem:range}, $\lambda>B_b-2=3b^2+3b+1\ge19$ for $|b|\ge3$ and for
$b=2$. For $b=-2$ and $b=-1$, the polynomial \eqref{eq:etamin} is
$X^3-9X^2-3X-1$ and $X^3-3X^2-1$ respectively. Each has one real root and is
negative at $X=9.3$ and at $X=3.1$ respectively, so $\lambda>9.3$ and
$\lambda>3.1$.

Hence $\lambda^m>3.1^{31}>10^{15}$ in all cases. By Lemma~\ref{iso:identity},
$|\cos2m\psi-\frac12|=|R_m|/2\le3\lambda^{-m/2}+5\lambda^{-m}\le
3.5\lambda^{-m/2}$. Lemma~\ref{iso:angle} first gives
$\varepsilon:=\operatorname{dist}(2m\psi,\mathcal{A})\le3.5\pi\lambda^{-m/2}<\frac1{10}$,
and then $\varepsilon\le(3.5/0.815)\lambda^{-m/2}<4.3\lambda^{-m/2}$.

As $2m\psi\equiv\pm m\vartheta\pmod{2\pi}$ and $\mathcal{A}=-\mathcal{A}$, there is
$k\in\intZ$ with $|m\vartheta\mp\frac\pi3-2\pi k|\le4.3\lambda^{-m/2}$ for one of the two
signs. Multiply by $3$ and put $j=6k\pm1$, which is odd. Then $|j\pi-3m\vartheta|\le12.9\lambda^{-m/2}<1$. Since
$0<3m\vartheta<3m\pi$, this forces $-1<j<3m+1$, so $1\le j\le3m$.
\end{proof}

We apply a theorem of Laurent, Mignotte and Nesterenko, which is
Th\'eor\`eme~3 of \cite{LMN1995}. Let $\alpha$ be an algebraic number of absolute
value $1$ that is not a root of unity, let $b_1,b_2$ be positive integers, and let
$\mathcal{L}=b_1\mathrm{i}\pi-b_2\log\alpha$. Put
\begin{gather*}
  d=[\ratQ(\alpha):\ratQ]/2,\qquad a=\max\{20,\ 10.98\,|\log\alpha|+d\,h(\alpha)\},\\
  H=\max\Bigl\{17,\ \frac{\sqrt d}{10},\
     d\log\Bigl(\frac{b_1}{2a}+\frac{b_2}{68.9}\Bigr)+2.35\,d+5.03\Bigr\}.
\end{gather*}
Then $\log|\mathcal{L}|\ge-8.87\,aH^2$. Here $h$ is the absolute logarithmic height,
the logarithm may be any determination \cite[\S\S1--2]{LMN1995}, and $d$ is
written $D$ in \cite{LMN1995}; $a$ and $H$ are named as there.

\begin{proposition}\label{iso:lmn}
Under the hypotheses of Lemma~\ref{iso:linear}, $G(m,\log\lambda)\le0$, where
\begin{gather*}
  G(m,x)=\frac m2\,x-\log13-8.87\bigl(10.98\pi+\tfrac32x\bigr)\bar H(m)^2,\\
  \bar H(m)=\max\Bigl\{17,\ 3\log\Bigl(\frac{3m}{40}+\frac{3m}{68.9}\Bigr)+12.08\Bigr\}.
\end{gather*}
\end{proposition}

\begin{proof}
Let $\alpha\in\{\gamma,\bar\gamma\}$ be the one with argument $\vartheta$, and take
$\log\alpha=\mathrm{i}\vartheta$. By Lemma~\ref{iso:beta}, $\alpha$ is not a root of
unity, $d\le3$ and $h(\alpha)=\frac12\log\lambda$. With $b_1=j$ and $b_2=3m$ from
Lemma~\ref{iso:linear}, we get $\mathcal{L}=\mathrm{i}(j\pi-3m\vartheta)$, which is non-zero
as $\alpha$ is not a root of unity. Then $20\le a\le10.98\pi+\frac32\log\lambda$, since
$10.98\pi>20$. The expression in $H$ increases with $d$, because
$b_2/68.9\ge93/68.9>1$. Using $b_1\le3m$, $a\ge20$ and $d\le3$, we get $H\le\bar H(m)$.
Combining \eqref{eq:linform} with the theorem,
$\log13-\frac m2\log\lambda\ge\log|\mathcal{L}|\ge-8.87\,aH^2\ge
-8.87(10.98\pi+\frac32\log\lambda)\bar H(m)^2$.
\end{proof}

The parameter $\log\lambda$ does not cancel in $G$, but it enters $G$ with the
positive coefficient $\frac m2-13.305\,\bar H(m)^2$ once $m$ is large. Hence
Proposition~\ref{iso:lmn} bounds $m$ independently of $b$. To bound $b$ as well,
we use the fact that $\eta_b$ itself is of family type, so $2\psi$ is close to
$\mathcal{A}$, but not too close.

\begin{lemma}\label{iso:eps}
Let $|b|\ge10^4$ and $\varepsilon_b=\operatorname{dist}(2\psi,\mathcal{A})$. Then
\[
  (2.4-5\cdot10^{-4})\,\lambda^{-1/2}\ \le\ \varepsilon_b\ \le\ \frac{2.13}{|b|} .
\]
If moreover $m\ge31$ and $\eta_b^{\,m}$ is of family type, then $m\ge\nu|b|$,
where $\nu=(\frac\pi3-10^{-9})/2.13$.
\end{lemma}

\begin{proof}
Here $\lambda>3b^2+3b+1\ge3b^2-3|b|+1>10^8$ by Lemma~\ref{lem:range}, so
$\lambda^{-1/2}<10^{-4}$.
Lemma~\ref{iso:identity} with $m=1$ applies, since $\eta_b$ is of family type. It
gives $r:=|\cos2\psi-\frac12|=|R_1|/2\le3\lambda^{-1/2}+5\lambda^{-1}$.

Upper bound. Lemma~\ref{iso:angle} first gives
$\varepsilon_b\le\pi r\le\pi(3+5\cdot10^{-4})10^{-4}<\frac1{10}$, and then
$\varepsilon_b\le r/0.815\le(3/0.815)(1+\frac53\cdot10^{-4})\lambda^{-1/2}$.
Moreover $\lambda\ge3b^2(1-|b|^{-1})\ge3(1-10^{-4})b^2$. Together these give
$\varepsilon_b|b|\le(3/0.815)(1+\frac53\cdot10^{-4})/\sqrt{3(1-10^{-4})}<2.13$.

Lower bound. Since $0<2\psi<2\pi$ lies within $\frac1{10}$ of $\frac\pi3$ or
$\frac{5\pi}3$, $\psi$ lies within $\frac1{20}$ of $\frac\pi6$ or $\frac{5\pi}6$.
Hence $|\cos\psi|\ge\frac{\sqrt3}2\cos\frac1{20}-\frac12\sin\frac1{20}>0.8$. So
$r\ge3\lambda^{-1/2}|\cos\psi|-5\lambda^{-1}\ge2.4\lambda^{-1/2}-5\lambda^{-1}$, and
$\varepsilon_b\ge r$ by Lemma~\ref{iso:angle}.

Finally, write $2\psi=x_0+y$ with $x_0\in\mathcal{A}$ and $|y|=\varepsilon_b$.
Lemma~\ref{iso:linear} gives $x_1\in\mathcal{A}$ with
$|2m\psi-x_1|\le t_m:=4.3\lambda^{-m/2}$. Since $mx_0$ and $x_1$ lie in
$\frac\pi3\intZ$, the distance from $my$ to $\frac\pi3\intZ$ is at most $t_m$.
Now $t_m\le4.3\lambda^{-15.5}<\varepsilon_b\le m|y|$, so the nearest point
of $\frac\pi3\intZ$ to $my$ is not $0$. Hence
$m\varepsilon_b\ge\frac\pi3-t_m>\frac\pi3-10^{-9}$, and the upper bound for
$\varepsilon_b$ gives the claim.
\end{proof}

\begin{proposition}\label{iso:regions}
Let $m_{\mathrm I}=630000$, $m_{\mathrm{III}}=1310000$ and $N=163000$. Then $G(m,x)>0$ in each of the
regions
\begin{enumerate}
\item $m\ge m_{\mathrm I}$ and $x\ge\log9.3$;
\item $m\ge\nu N$ and $x\ge\log(3N^2-3N+1)$;
\item $m\ge m_{\mathrm{III}}$ and $x=\log\eta_{-1}=1.1326\ldots$.
\end{enumerate}
\end{proposition}

\begin{proof}
For $m\ge50$, $\bar H(m)=3\log(cm)+12.08$ with $c=\frac3{40}+\frac3{68.9}$. So
$\bar H(m)>6$, and both $\bar H(m)^2/m$ and $\bar H(m)/m$ decrease. We have
$\partial G/\partial x=\frac m2-13.305\,\bar H(m)^2$. Also
$\partial G/\partial m=\frac x2-53.22\,(10.98\pi+\frac32x)\bar H(m)/m$, which is
positive when $m>106.44\,(10.98\pi/x+\frac32)\bar H(m)$, and the right-hand side
decreases in $x$. So if both partial derivatives are positive at the corner
$(m_*,x_*)$ of a region, they are positive on the whole region. Then
$G(m,x)\ge G(m,x_*)\ge G(m_*,x_*)$. It remains to check, at each of the three
corners $(m_{\mathrm I},\log9.3)$, $(\nu N,\log(3N^2-3N+1))$ and
$(m_{\mathrm{III}},x_-)$, where $x_-$ is a rational lower bound for
$\log\eta_{-1}$, that both partial derivatives are positive and that $G>0$. We
do this with rational enclosures of $\pi$ and of the logarithms, which give
$G(m_{\mathrm I},\log9.3)>156$, $G(\nu N,\log(3N^2-3N+1))>4430$ and $G(m_{\mathrm{III}},x_-)>4073$.
The first margin is narrow: $156$ is small against terms of size about $7\cdot10^5$,
and the least $m$ with $G(m,\log9.3)>0$ is about $629839$.
\end{proof}

It remains to treat the finite range $|b|<N$, $31\le m<m_{\mathrm I}$, with $m_{\mathrm{III}}$ in place of $m_{\mathrm I}$
for $b=-1$. Put $\xi=3\vartheta/\pi$ and
write $\|t\|$ for the distance from $t\in\mathbb{R}$ to $\intZ$. By
\eqref{eq:linform}, a power $\eta_b^{\,m}$ of family type with $m\ge31$ satisfies
\begin{equation}\label{eq:xi}
  \|m\xi\|\le\frac{13}{\pi}\lambda^{-m/2} .
\end{equation}

\begin{proposition}\label{iso:cert}
Let $b$ be an integer with $2\le|b|<N$. Then
$\|m\xi\|>\frac{13}\pi\lambda^{-15}$ for every $m$ with $31\le m\le m_{\mathrm I}$. For
$b=-1$, $\|m\xi\|>\frac{13}\pi\lambda^{-15}$ for $31\le m\le1000$, and
$\|m\xi\|>\frac{13}\pi\lambda^{-500}$ for $1000\le m\le m_{\mathrm{III}}$.
\end{proposition}

By \eqref{eq:xi}, this excludes every $m$ in these ranges, because
$\lambda^{-m/2}\le\lambda^{-15}$ for $m\ge31$ and $\lambda^{-m/2}\le\lambda^{-500}$
for $m\ge1000$.

\begin{proof}
For each of the $325997$ values of $b$ we verify the following in exact rational
arithmetic.
\begin{enumerate}
\item Choose rationals $\lambda_-<\lambda_+$ at which \eqref{eq:etamin} takes
  opposite signs, so that $\lambda\in[\lambda_-,\lambda_+]$. The width is $2^{-317}$.
\item Bound $\cos\vartheta$ on this interval by \eqref{eq:cosomega}. Evaluate the
  right-hand side at the endpoints and add a term $(w^2/8)\max|f''|$, where
  $f(y)=\frac12y(B_b-y)^2-1$ and $w=\lambda_+-\lambda_-$.
\item Choose rationals $0<t_-<t_+<\pi$ with $\cos t_+$ and $\cos t_-$ on either side
  of the enclosure of $\cos\vartheta$, so that $\vartheta\in[t_-,t_+]$. Enclose
  $\cos$ by its Taylor series, with the remainder bounded, and $\pi$ by
  Machin's formula. Put $\xi_-=3t_-/\pi_+$ and $\xi_+=3t_+/\pi_-$, where
  $\pi_-<\pi<\pi_+$.
\item Let $q_k\le m_{\mathrm I}<q_{k+1}$ be consecutive denominators of convergents of the
  rational number $\xi_-$. By the best-approximation property of convergents
  \cite[Theorems~16 and~17]{Khinchin1964}, $\|m\xi_-\|\ge\|q_k\xi_-\|$ for $1\le m<q_{k+1}$:
  the least $m$ minimizing $\|m\xi_-\|$ over this range is the denominator of a best
  approximation of the second kind, hence of a convergent $p_j/q_j$ with $j\le k$
  (Theorem~16), and $|q_j\xi_--p_j|\ge|q_k\xi_--p_k|$ by Theorem~17, applied to
  $p_k/q_k$. Theorem~17 has one exception, $\xi_-\in\intZ+\frac12$, which
  concerns $k=0$; then $j=k$, and the inequality holds trivially. Hence
  $\|m\xi\|\ge\|q_k\xi_-\|-m_{\mathrm I}(\xi_+-\xi_-)$ for $m\le m_{\mathrm I}$. Check that
  the right-hand side exceeds $\frac{13}{\pi_-}\lambda_-^{-15}$.
\end{enumerate}
For $b=-1$, we carry out step~(iv) twice: with $1000$ and exponent $15$, and
with $m_{\mathrm{III}}$ and exponent $500$. A single step with $m_{\mathrm{III}}$ and exponent $15$ does not
suffice. It gives $\|q_k\xi_-\|=1.6\cdot10^{-7}$, while
$\frac{13}\pi\lambda^{-15}=1.7\cdot10^{-7}$.
\end{proof}

\begin{proof}[Proof of Theorem~\ref{iso:main}]
For $2\le m\le30$ use Lemma~\ref{iso:small}; the pair $(-1,4)$ does occur, by
Example~\ref{iso:ex}. Let $m\ge31$, and suppose that $\eta_b^{\,m}$ is of family
type. By Proposition~\ref{iso:lmn}, $G(m,\log\lambda)\le0$. If $|b|\ge N$, then
$m\ge\nu|b|\ge\nu N$ by Lemma~\ref{iso:eps}, and
$\log\lambda\ge\log(3b^2-3|b|+1)\ge\log(3N^2-3N+1)$. This contradicts
Proposition~\ref{iso:regions}(ii). If $2\le|b|<N$, then $\lambda>9.3$ by
Lemma~\ref{iso:linear}, so Proposition~\ref{iso:regions}(i) gives $m<m_{\mathrm I}$. If
$b=-1$, Proposition~\ref{iso:regions}(iii) gives $m<m_{\mathrm{III}}$. In either case
Proposition~\ref{iso:cert} contradicts~\eqref{eq:xi}.
\end{proof}

\begin{remark}\label{iso:status}
Lemmas~\ref{iso:small} and~\ref{iso:six} and Propositions~\ref{iso:regions}
and~\ref{iso:cert} are computations. They are exact: integer and
rational arithmetic throughout. We use floating point only to propose a candidate,
such as a root or an angle, which we then verify by a sign change or an
interval inequality. Everything else in the proof is written out above. The
computation of Proposition~\ref{iso:cert} takes about three minutes.
\end{remark}

\subsection{Consequences}\label{sec:isocons}

\begin{corollary}\label{iso:cor}
Let $b\neq c$ be integers different from $0,1$, with $K_c=K_b$. If
$n(b)\mid n(c)$, then $(b,c)=(-1,5)$. In particular:
\begin{enumerate}
\item if $n(b)=1$ and $b\neq-1$, then $K_b$ is isomorphic to no $K_c$ with
  $c\neq b$;
\item $K_{-1}=K_5$, and this field is $K_c$ for no other $c$.
\end{enumerate}
\end{corollary}

\begin{proof}
By Lemma~\ref{iso:reform}, $\eta_c=\eta_0^{\,n(c)}$ with $n(c)\neq n(b)$. If
$n(b)\mid n(c)$, then $\eta_c=\eta_b^{\,m}$ with $m=n(c)/n(b)\ge2$, a power of
family type. Theorem~\ref{iso:main} gives $(b,m)=(-1,4)$, and then
$\eta_c=\eta_{-1}^{\,4}=\eta_5$ (Example~\ref{iso:ex}), so $c=5$. Part~(i)
follows, since $n(b)=1$ divides every $n(c)$. So does part~(ii), because
$n(-1)=1$.
\end{proof}

By \cite{KimuraUmemoto}, $n(b)=1$ holds in particular when $3\nmid b$, $b^3-4$
is squarefree and either $b\ge7$ or $b\le-1$, a set of positive density.
Within that set the non-isomorphism follows from Lemma~\ref{iso:reform} alone,
since $n(b)=n(b')=1$ for any two members $b,b'$. The
point of Corollary~\ref{iso:cor} is that the other members of the family are
excluded too, except $K_5=K_{-1}$ for $b=-1$. Using Theorems~\ref{thmB}, \ref{thmC} and~\ref{thmD}, we show in
Corollary~\ref{iso:far} that any other coincidence involves a parameter outside
the range of Theorem~\ref{thmC}. We need one more computation.

\begin{lemma}\label{iso:six}
Let $b\in\{-75,-3,3,5,9,96\}$, the $b$ with $n(b)\ge3$ in Theorem~\ref{thmC}, and
let $\eta_0$ be the fundamental unit of $K_b$. For $1\le j\le30$, $\eta_0^{\,j}$ is
of family type only for $j=n(b)$ and, when $b=5$, also for $j=1$.
\end{lemma}

\begin{proof}
We compute $\eta_0$ in $K_b$ as the $n(b)$-th root of $\eta_b$ that is greater
than $1$, and single it out exactly, as the root of norm $+1$: for a unit $u$ of
$K_b$, the real embedding has the sign of $\Norm(u)=u\,|\sigma(u)|^2$. We then
test the characteristic polynomials of its powers against \eqref{eq:D}.
\end{proof}

\begin{corollary}\label{iso:far}
Let $b\neq b'$ be integers different from $0,1$, with $K_b=K_{b'}$ and
$\{b,b'\}\neq\{-1,5\}$. Then there are a $c\in\{b,b'\}$ with $|c|>10^{30}$ and a
prime $q\ge5$ dividing $n(c)$.
\end{corollary}

\begin{proof}
By Corollary~\ref{iso:cor}, applied both ways, neither of $n(b),n(b')$ divides
the other. In particular both are at least $2$.

First let $|b|\le10^{30}$ and $|b'|\le10^{30}$. Then, by Theorem~\ref{thmC},
$n(b)$ and $n(b')$ lie in $\{2,4,6,8,13\}$ and are distinct. So one of them, say $n(b)$, is at least
$3$. Then $b$ is one of the six parameters of Lemma~\ref{iso:six}, and
$\eta_0^{\,n(b')}$ is of family type with $2\le n(b')\le13$ and $n(b')\neq n(b)$.
The lemma excludes this.

So we may assume $|b'|>10^{30}$. If a prime $q\ge5$ divides $n(b')$, take $c=b'$.
Otherwise $n(b')=2^s3^t$. Here $t\ge1$ would give $b'=9$ by
Theorem~\ref{thmB}(ii), and $s\ge2$ would give $b'\in\{-75,-3,5,96\}$ by
Theorem~\ref{thmD}. Hence $n(b')=2$. As $n(b')\nmid n(b)$, $n(b)$ is odd and at
least $3$. It is prime to $3$, since otherwise $b=9$, and $n(9)=6$ is even. So a
prime $q\ge5$ divides $n(b)$. Finally, if $|b|\le10^{30}$, then $b=3$ and
$n(b)=13$ by Theorem~\ref{thmC}, and $\eta_0^{\,2}$ would be of family type in
$K_3$, which Lemma~\ref{iso:six} excludes. Hence $|b|>10^{30}$, and we take
$c=b$.
\end{proof}

\begin{remark}
Corollary~\ref{iso:far} inherits the status of Theorem~\ref{thmD}
(Remark~\ref{rem:status}). Corollary~\ref{iso:cor} does not use Theorems~\ref{thmC}
or~\ref{thmD}. By Corollary~\ref{iso:far}, a coincidence $K_b=K_{b'}$ other
than $K_{-1}=K_5$, if one exists, requires $q\mid n(c)$ for some
$c\in\{b,b'\}$ and a prime $q\ge5$, a case to which the curves $C_q$ of genus
$\binom{q-1}2\ge6$ are relevant. We do not know whether there is any coincidence
besides $K_{-1}=K_5$. By Corollary~\ref{iso:finite} below, there are at most
finitely many. Proposition~\ref{iso:square} below restates the question
as a Diophantine one.
\end{remark}

The proofs of Lemmas~\ref{iso:identity} and~\ref{iso:beta} use only that
$\lambda$ is a unit of $K_b$ greater than $1$. So they hold for every such unit
$\lambda$, with $\mu$, $\psi$, $\gamma$ and $\vartheta$ defined from $\lambda$ as
above and with $\eta_b^{\,m}$ replaced by $\lambda^m$. Since
Lemma~\ref{lem:smyth} bounds $\lambda$ from below independently of $b$, the
argument of \S\ref{sec:isopow} then bounds the exponent $m$ independently of
$b$.

\begin{proposition}\label{iso:uniform}
Let $b\neq0,1$, let $\lambda>1$ be a unit of $K_b$, and let $m\ge1$. If
$\lambda^m$ is of family type, then $m<M_0$, where $M_0=6082138$.
\end{proposition}

\begin{proof}
We may assume $m\ge124$. By Lemma~\ref{lem:smyth}, $\lambda\ge\rho>1.3247$, so
$\lambda^m>1.3247^{124}>10^{15}$. The proof of Lemma~\ref{iso:linear} uses only
$\lambda^m>10^{15}$ from the point where this bound appears. So it gives an odd
$j$ with $1\le j\le3m$ and $|j\pi-3m\vartheta|\le13\lambda^{-m/2}$. The proof of
Proposition~\ref{iso:lmn} then applies without change, since
$b_2=3m\ge372>68.9$, and gives $G(m,\log\lambda)\le0$.

On the other hand, $G(m,x)>0$ for $m\ge M_0$ and $x\ge\log1.3247$. We prove this
as in the proof of Proposition~\ref{iso:regions}, at the corner $(M_0,x_0)$, where
$x_0$ is a rational lower bound for $\log1.3247$. Both partial derivatives are
positive there, and $G(M_0,x_0)>0.005$. As for Proposition~\ref{iso:regions},
we check this with rational enclosures of $\pi$ and of the logarithms
(Remark~\ref{iso:status}). Since $\log\lambda>\log1.3247$, we get
$m<M_0$.
\end{proof}

The margin is narrow: $G(M_0-1,x)<0$ for $x=\log1.3247$, so this argument gives
no smaller bound.

\begin{corollary}\label{iso:finite}
\begin{enumerate}
\item $n(b)<M_0$ for every $b\neq0,1$.
\item Only finitely many integers $b\neq0,1$ have $n(b)\notin\{1,2\}$.
\item Only finitely many pairs $\{b,b'\}$ of distinct integers different from
  $0,1$ satisfy $K_b\cong K_{b'}$.
\end{enumerate}
\end{corollary}

\begin{proof}
(i) Let $\eta_0$ be the fundamental unit of $K_b$. Then $\eta_0^{\,n(b)}=\eta_b$ is
of family type, so Proposition~\ref{iso:uniform} applies with $\lambda=\eta_0$ and
$m=n(b)$.

(ii) If $n(b)\ge3$, then $n(b)$ is divisible by $3$, by $4$ or by a prime
$q\ge5$, and $q<M_0$ by (i). By Theorem~\ref{thmB}(ii), only $b=9$ has
$3\mid n(b)$. By Theorem~\ref{thmB}(iii), each $q\ge4$ divides $n(b)$ for only
finitely many $b$. So the set $S$ of $b$ with $n(b)\ge3$ is a finite union of
finite sets.

(iii) Let $K_b\cong K_{b'}$ with $b\neq b'$ and $\{b,b'\}\neq\{-1,5\}$. By
Lemma~\ref{iso:reform}, $K_b=K_{b'}$ and $n(b)\neq n(b')$. By
Corollary~\ref{iso:cor}, applied both ways, neither of $n(b),n(b')$ divides the
other. So both are at least $2$, and as they are distinct, one of them, say
$n(b)$, is at least $3$. Thus $b\in S$. By Proposition~\ref{iso:bij}, the $c$
with $K_c=K_b$ correspond to the $j\ge1$ for which $\eta_0^{\,j}$ is of family
type, and $j<M_0$ by Proposition~\ref{iso:uniform}. So, for each $b\in S$, fewer
than $M_0$ integers $b'$ are possible.
\end{proof}

Parts (i) and (ii) of Corollary~\ref{iso:finite} are Theorem~\ref{thmB}(iv), and
part (iii) is the last assertion of Theorem~\ref{thmE}(iii). Part (i) is
effective. Parts (ii) and (iii) are not, since they rest on
Theorem~\ref{thmB}(iii). To determine the set $S$ by this route, one would need
the integral points of $C_q$ for every prime $5\le q<M_0$.

Finally, the coincidences $K_b=K_{b'}$ can be read off from $b$ and $b'$ alone.

\begin{proposition}\label{iso:square}
Let $b$ and $b'$ be integers different from $0$ and $1$. Then $K_b=K_{b'}$ if
and only if $b(b^3-4)\,b'(b'^3-4)$ is a square.
\end{proposition}

\begin{proof}
If $b=b'$, there is nothing to prove, so let $b\neq b'$. Put $n=b^3-2$ and
$\zeta=(\theta^2-2b)/b$, where $\theta$ is the real root of $f_b$. From
$\theta^3=3b\theta+b^3$ we get $\zeta^3-3\zeta=n$. If $\zeta$ were rational,
then $\theta^2=b\zeta+2b$ would contradict $[\ratQ(\theta):\ratQ]=3$. So
$K_b=\ratQ(\zeta)$, and $g_n(X)=X^3-3X-n$ is the minimal polynomial of $\zeta$.
Moreover, $|n|\ge3$ and $n^2-4=b^2\cdot b(b^3-4)$. The same holds for $b'$,
with $n'=b'^3-2$.

Suppose that $K_b=K_{b'}$. The discriminant of $g_n$ is $-27(n^2-4)$. It
differs from the discriminant of $K_b$ by a nonzero square factor, and so does
the discriminant of $g_{n'}$. So $(n^2-4)(n'^2-4)$ is a square, and hence so is
$b(b^3-4)\,b'(b'^3-4)$.

Conversely, suppose that $b(b^3-4)\,b'(b'^3-4)$ is a square. Let $D$ be the
squarefree part of $b(b^3-4)$. It is positive, and $D\neq1$, since
$n^2-4=s^2$ would force $(n-s)(n+s)=4$ and $|n|=2$. Write $n^2-4=s^2D$ and
$n'^2-4=s'^2D$ with positive integers $s,s'$, and let $k=\ratQ(\sqrt D)$, a real
quadratic field. Then $\beta=(n+s\sqrt D)/2$ and $\beta'=(n'+s'\sqrt D)/2$ are
units of $k$ of norm $1$. Neither is a cube in $k$. Indeed, if $\beta=\gamma^3$
with $\gamma\in k$, then $\gamma$ has norm $1$, and
$(\gamma+\gamma^{-1})^3-3(\gamma+\gamma^{-1})=\gamma^3+\gamma^{-3}=n$. So the
rational number $\gamma+\gamma^{-1}$ would be a root of $g_n$.

The units of norm $1$ of $k$ are $\pm\varepsilon^j$ with $j\in\intZ$, for a
unit $\varepsilon>1$. So $\beta=\pm\varepsilon^t$ and
$\beta'=\pm\varepsilon^{t'}$, with $3\nmid t$ and $3\nmid t'$. Choose
$e\in\{1,-1\}$ with $t\equiv et'\pmod3$. Since $-1=(-1)^3$, we get
$\beta=x^3\beta'^{\,e}$ with $x\in k$. Here $\beta'^{-1}$ is the conjugate of
$\beta'$. By a result of Komatsu \cite{Komatsu2001}, in the form stated by Kishi
\cite[Proposition~6.3]{Kishi2000}, $g_n$ and $g_{n'}$ therefore have the same
splitting field $M$. The discriminant of $g_n$
is negative, so $\mathrm{Gal}(M/\ratQ)\cong S_3$. The cubic subfields of $M$
are conjugate, so $K_b\cong K_{b'}$, and $K_b=K_{b'}$ by
Lemma~\ref{iso:reform}.
\end{proof}

For example, $b(b^3-4)$ is $5$ for $b=-1$ and $605=5\cdot11^2$ for $b=5$. So
Proposition~\ref{iso:square} gives $K_{-1}=K_5$ again. The discriminant of
$K_b$ lies in $-3b(b^3-4)\,\ratQ^{\times2}$. So Proposition~\ref{iso:square}
says that $K_b$ is determined by the quadratic subfield of its Galois closure.
This fails for cubic fields in general. For instance, $\ratQ(\sqrt[3]{2})$ and
$\ratQ(\sqrt[3]{3})$ are not isomorphic, but both Galois closures contain
$\ratQ(\sqrt{-3})$. In the proof, the Galois closure of $K_b$ comes from a cube
root of a unit of the real quadratic field $\ratQ(\sqrt{b(b^3-4)})$. The units
of norm $1$ of this field give only one cubic field.

Proposition~\ref{iso:square} turns the question of further coincidences into a
Diophantine one. Consider the equation
\[
  w^2=b(b^3-4)\,b'(b'^3-4),\qquad b\neq b',\quad b,b'\notin\{0,1\}.
\]
By Corollary~\ref{iso:finite}(iii), it has only finitely many solutions in
integers. By Corollary~\ref{iso:far}, its only solutions with
$|b|,|b'|\le10^{30}$ have $\{b,b'\}=\{-1,5\}$ and $w=\pm55$. The first
statement is not effective. The second inherits the status of
Theorem~\ref{thmD}.

\section*{The computations}

We did the sweep of \S\ref{sec:sweep}, the $2$-descent of \S\ref{sec:c4curve}
and the sieve of \S\ref{sec:sieve} in PARI/GP \cite{PARI2}. The descent uses
the class and unit groups of $K$ and of the sextic fields cut out by the three
quadratic factors of \eqref{eq:X}. Two of the factors cut out the same field,
so there are two such fields. All three class and unit groups are certified, so
the descent does not assume the generalized Riemann hypothesis. We did the
$p$-adic computation of \S\ref{sec:qc} in SageMath \cite{SageMath}, with Coleman
integrals computed by our own implementation of Tuitman's algorithm
\cite{Tuitman2016,BalakrishnanTuitman2020}. We computed both regular models in
SageMath, the one at $\ell=3$ from \cite{Dokchitser2021} and the one at $\ell=2$
from \cite{Liu2026}; the latter agrees with what Magma \cite{Magma} returns. We
did the elliptic-logarithm computation of \S\ref{sec:elllog} in PARI/GP, in exact
arithmetic. We
did the computations of \S\ref{sec:iso} in PARI/GP, in exact integer and
rational arithmetic. Every
script reports a pass-or-fail verdict for each
assertion it is meant to establish; the scripts are available in a public
repository \cite{KimuraCode2}, which records for each script the statements
it checks.

\section*{Use of artificial intelligence tools}
During the preparation of this work the authors used Anthropic's Claude in order
to generate the PARI/GP and SageMath scripts on which \S\ref{sec:sweep},
\S\ref{sec:c4} and \S\ref{sec:iso} rest, including the implementation of
Tuitman's algorithm used in \S\ref{sec:c4}, to cross-check the output of the numerical experiments, to assist with
literature searches, and to improve the English of the manuscript. The
mathematical content, including every proof, is the authors' own. The authors
reviewed and verified all AI-generated output, including every reference and
all final results, edited the content as needed, and take full responsibility
for the content of this article.

\section*{Statements and Declarations}

\subsection*{Competing interests}
The authors have no competing interests to declare that are relevant to the
content of this article.

\subsection*{Funding}
No funding was received for conducting this study.

\subsection*{Data availability}
No datasets were generated or analysed during the current study beyond the
output of the computations described in the section ``The computations'', which
is reproduced by the scripts listed under Code availability.

\subsection*{Code availability}
The PARI/GP and SageMath scripts that carry out the computations of this
article, each reporting a pass-or-fail verdict for the statements it checks,
are available in a public repository \cite{KimuraCode2}.

\bibliography{../biblio}

\providecommand{\bysame}{\leavevmode\hbox to3em{\hrulefill}\thinspace}
\providecommand{\doi}[1]{\url{https://doi.org/#1}}
\providecommand{\MR}{\relax\ifhmode\unskip\space\fi MR }
\providecommand{\MRhref}[2]{%
  \href{http://www.ams.org/mathscinet-getitem?mr=#1}{#2}
}
\providecommand{\href}[2]{#2}
\begin{thebibliography}{10}

\bibitem{BakerPayneRabinoff2013}
Matthew Baker, Sam Payne, and Joseph Rabinoff, \emph{On the structure of
  non-{A}rchimedean analytic curves}, Tropical and non-{A}rchimedean geometry,
  Contemp. Math., vol. 605, Amer. Math. Soc., Providence, RI, 2013,
  pp.~93--121. \doi{10.1090/conm/605/12113}

\bibitem{BalakrishnanDogra2018}
Jennifer~S. Balakrishnan and Netan Dogra, \emph{Quadratic {C}habauty and
  rational points {I}: $p$-adic heights}, Duke Math. J. \textbf{167} (2018),
  no.~11, 1981--2038, Preprint: arXiv:1601.00388.
  \doi{10.1215/00127094-2018-0013}

\bibitem{BalakrishnanDogra2021}
\bysame, \emph{Quadratic {C}habauty and rational points {II}: generalised
  height functions on {S}elmer varieties}, Int. Math. Res. Not. IMRN (2021),
  no.~15, 11923--12008, Preprint: arXiv:1705.00401. \doi{10.1093/imrn/rnz362}

\bibitem{BDMTV2019}
Jennifer~S. Balakrishnan, Netan Dogra, J.~Steffen M\"uller, Jan Tuitman, and
  Jan Vonk, \emph{Explicit {C}habauty--{K}im for the split {C}artan modular
  curve of level 13}, Ann. of Math. (2) \textbf{189} (2019), no.~3, 885--944,
  Preprint: arXiv:1711.05846. \doi{10.4007/annals.2019.189.3.6}

\bibitem{BDMTV2023}
\bysame, \emph{Quadratic {C}habauty for modular curves: algorithms and
  examples}, Compos. Math. \textbf{159} (2023), no.~6, 1111--1152, Preprint:
  arXiv:2101.01862. \doi{10.1112/S0010437X23007170}

\bibitem{BalakrishnanTuitman2020}
Jennifer~S. Balakrishnan and Jan Tuitman, \emph{Explicit {C}oleman integration
  for curves}, Math. Comp. \textbf{89} (2020), no.~326, 2965--2984, Preprint:
  arXiv:1710.01673. \doi{10.1090/mcom/3542}

\bibitem{Magma}
Wieb Bosma, John Cannon, and Catherine Playoust, \emph{The {M}agma algebra
  system. {I}. {T}he user language}, J. Symbolic Comput. \textbf{24} (1997),
  no.~3-4, 235--265. \doi{10.1006/jsco.1996.0125}

\bibitem{BruinStoll2010}
Nils Bruin and Michael Stoll, \emph{The {M}ordell--{W}eil sieve: proving
  non-existence of rational points on curves}, LMS J. Comput. Math. \textbf{13}
  (2010), 272--306. \doi{10.1112/S1461157009000187}

\bibitem{David1995}
Sinnou David, \emph{Minorations de formes lin\'eaires de logarithmes
  elliptiques}, M\'em. Soc. Math. France (N.S.) \textbf{62} (1995), 143 pp.
  \doi{10.24033/msmf.376}

\bibitem{Dokchitser2021}
Tim Dokchitser, \emph{Models of curves over discrete valuation rings}, Duke
  Math. J. \textbf{170} (2021), no.~11, 2519--2574.
  \doi{10.1215/00127094-2020-0079}

\bibitem{MR718935}
G.~Faltings, \emph{Endlichkeitss\"atze f\"ur abelsche {Variet\"aten} \"uber
  {Zahlk\"orpern}}, Invent. Math. \textbf{73} (1983), no.~3, 349--366, Erratum:
  \emph{Invent. Math.} \textbf{75} (1984), 381. \doi{10.1007/BF01388432}

\bibitem{GebelPethoZimmer1994}
J.~Gebel, A.~Peth\H{o}, and H.~G. Zimmer, \emph{Computing integral points on
  elliptic curves}, Acta Arith. \textbf{68} (1994), no.~2, 171--192.
  \doi{10.4064/aa-68-2-171-192}

\bibitem{Hoshi2011}
Akinari Hoshi, \emph{On correspondence between solutions of a family of cubic
  {T}hue equations and isomorphism classes of the simplest cubic fields}, J.
  Number Theory \textbf{131} (2011), no.~11, 2135--2150.
  \doi{10.1016/j.jnt.2011.05.001}

\bibitem{HoshiMiyake2010}
Akinari Hoshi and Katsuya Miyake, \emph{A note on the field isomorphism problem
  of {$X^3+sX+s$} and related cubic {T}hue equations}, Interdiscip. Inf. Sci.
  \textbf{16} (2010), no.~1, 45--54. \doi{10.4036/iis.2010.45}

\bibitem{Khinchin1964}
A.~Ya. Khinchin, \emph{Continued fractions}, The University of Chicago Press,
  Chicago and London, 1964.

\bibitem{Khovanskii1977}
A.~G. Khovanski\u{\i}, \emph{Newton polyhedra and toroidal varieties}, Funct.
  Anal. Appl. \textbf{11} (1977), 289--296. \doi{10.1007/BF01077143}

\bibitem{Khovanskii1978}
\bysame, \emph{Newton polyhedra and the genus of complete intersections},
  Funct. Anal. Appl. \textbf{12} (1978), no.~1, 38--46.
  \doi{10.1007/BF01077562}

\bibitem{KimuraCode2}
Iwao Kimura, \emph{{PARI/GP} and {SageMath} scripts for ``{On explicit units in
  a family of cubic fields and the arithmetic of curves arising from powers of
  these units}''}, GitHub repository, 2026, Version v1.3. Available from
  \url{https://github.com/iwaokimura/cubic-fields-unit-powers}.

\bibitem{KimuraUmemoto}
Iwao Kimura and Hikaru Umemoto, \emph{On an infinite family of cubic fields
  with explicit fundamental units}, 2026, Preprint: arXiv:2609.37389.

\bibitem{Kishi2000}
Yasuhiro Kishi, \emph{A constructive approach to {S}piegelung relations between
  3-ranks of absolute ideal class groups and congruent ones modulo {$(3)^2$} in
  quadratic fields}, J. Number Theory \textbf{83} (2000), no.~1, 1--49.
  \doi{10.1006/jnth.1999.2498}

\bibitem{Komatsu2001}
Toru Komatsu, \emph{On unramified cyclic cubic extensions of real quadratic
  fields}, Japan. J. Math. (N.S.) \textbf{27} (2001), no.~2, 353--386.
  \doi{10.4099/math1924.27.353}

\bibitem{KuribayashiKomiya1979}
Akikazu Kuribayashi and Kaname Komiya, \emph{On {W}eierstrass points and
  automorphisms of curves of genus three}, Algebraic geometry ({P}roc. {S}ummer
  {M}eeting, {C}openhagen, 1978), Lecture Notes in Math., vol. 732, Springer,
  Berlin, 1979, pp.~253--299. \doi{10.1007/BFb0066649}

\bibitem{Lang2002}
Serge Lang, \emph{Algebra}, third ed., Graduate Texts in Mathematics, vol. 211,
  Springer-Verlag, New York, 2002. \doi{10.1007/978-1-4613-0041-0}

\bibitem{LMN1995}
Michel Laurent, Maurice Mignotte, and Yuri Nesterenko, \emph{Formes lin\'eaires
  en deux logarithmes et d\'eterminants d'interpolation}, J. Number Theory
  \textbf{55} (1995), no.~2, 285--321. \doi{10.1006/jnth.1995.1141}

\bibitem{LevesqueWaldschmidt2015}
Claude L\'evesque and Michel Waldschmidt, \emph{A family of {T}hue equations
  involving powers of units of the simplest cubic fields}, J. Th\'eor. Nombres
  Bordeaux \textbf{27} (2015), no.~2, 537--563, Preprint: arXiv:1505.06708.
  \doi{10.5802/jtnb.913}

\bibitem{Liu2026}
Qing Liu, \emph{Desingularization of double covers of regular surfaces}, 2026,
  Preprint: arXiv:2504.16808v3.

\bibitem{RitzenthalerRomagny2018}
Christophe Ritzenthaler and Matthieu Romagny, \emph{On the {P}rym variety of
  genus 3 covers of genus 1 curves}, \'Epijournal G\'eom. Alg\'ebrique
  \textbf{2} (2018), Art.~2, Preprint: arXiv:1612.07033.
  \doi{10.46298/epiga.2018.volume2.3663}

\bibitem{Silverman1990}
Joseph~H. Silverman, \emph{The difference between the {W}eil height and the
  canonical height on elliptic curves}, Math. Comp. \textbf{55} (1990),
  no.~192, 723--743. \doi{10.1090/S0025-5718-1990-1035944-5}

\bibitem{Smyth1971}
C.~J. Smyth, \emph{On the product of the conjugates outside the unit circle of
  an algebraic integer}, Bull. London Math. Soc. \textbf{3} (1971), 169--175.
  \doi{10.1112/blms/3.2.169}

\bibitem{Steinberg1965}
Robert Steinberg, \emph{Regular elements of semi-simple algebraic groups},
  Publ. Math. Inst. Hautes \'Etudes Sci. \textbf{25} (1965), 49--80.
  \doi{10.1007/BF02684397}

\bibitem{Stichtenoth2009}
Henning Stichtenoth, \emph{Algebraic function fields and codes}, second ed.,
  Graduate Texts in Mathematics, vol. 254, Springer-Verlag, Berlin, 2009.
  \doi{10.1007/978-3-540-76878-4}

\bibitem{Stoll2001}
Michael Stoll, \emph{Implementing 2-descent for {J}acobians of hyperelliptic
  curves}, Acta Arith. \textbf{98} (2001), no.~3, 245--277.
  \doi{10.4064/aa98-3-4}

\bibitem{StroekerTzanakis1994}
R.~J. Stroeker and N.~Tzanakis, \emph{Solving elliptic {D}iophantine equations
  by estimating linear forms in elliptic logarithms}, Acta Arith. \textbf{67}
  (1994), no.~2, 177--196. \doi{10.4064/aa-67-2-177-196}

\bibitem{LMFDB}
{The LMFDB Collaboration}, \emph{The {L}-functions and modular forms database},
  \url{https://www.lmfdb.org}, 2026, [Online; accessed 19 November 2026].

\bibitem{PARI2}
{The PARI Group}, \emph{{PARI/GP} version \texttt{2.17.4}}, 2026, Available
  from \url{https://pari.math.u-bordeaux.fr/}.

\bibitem{SageMath}
{The Sage Developers}, \emph{{S}age{M}ath, the {S}age {M}athematics {S}oftware
  {S}ystem (version 10.7)}, 2025, Available from
  \url{https://www.sagemath.org}.

\bibitem{Tuitman2016}
Jan Tuitman, \emph{Counting points on curves using a map to {$\mathbf{P}^1$}},
  Math. Comp. \textbf{85} (2016), no.~298, 961--981, Preprint: arXiv:1402.6758.
  \doi{10.1090/mcom/2996}

\end{thebibliography}

\end{document}